\documentclass[twoside,leqno,10pt, A4]{amsart}
\usepackage{amsfonts}
\usepackage{amsmath}
\usepackage{amscd}
\usepackage{amssymb}
\usepackage{amsthm}
\usepackage{amsrefs}
\usepackage{latexsym}
\usepackage{mathrsfs}
\usepackage{bbm}
\usepackage{amscd}
\usepackage{amssymb}
\usepackage{amsthm}
\usepackage{amsrefs}
\usepackage{latexsym}
\usepackage{mathrsfs}
\usepackage{bbm}
\usepackage{enumerate}
\usepackage{graphicx}
\usepackage{braket}
\usepackage{color}
\begin{document}

\newtheorem{theorem}[subsection]{Theorem}
\newtheorem{proposition}[subsection]{Proposition}
\newtheorem{lemma}[subsection]{Lemma}
\newtheorem{corollary}[subsection]{Corollary}
\newtheorem{conjecture}[subsection]{Conjecture}
\newtheorem{prop}[subsection]{Proposition}
\newtheorem{defin}[subsection]{Definition}

\numberwithin{equation}{section}
\newcommand{\mr}{\ensuremath{\mathbb R}}
\newcommand{\mc}{\ensuremath{\mathbb C}}
\newcommand{\dif}{\mathrm{d}}
\newcommand{\intz}{\mathbb{Z}}
\newcommand{\ratq}{\mathbb{Q}}
\newcommand{\natn}{\mathbb{N}}
\newcommand{\comc}{\mathbb{C}}
\newcommand{\rear}{\mathbb{R}}
\newcommand{\prip}{\mathbb{P}}
\newcommand{\uph}{\mathbb{H}}
\newcommand{\fief}{\mathbb{F}}
\newcommand{\majorarc}{\mathfrak{M}}
\newcommand{\minorarc}{\mathfrak{m}}
\newcommand{\sings}{\mathfrak{S}}
\newcommand{\fA}{\ensuremath{\mathfrak A}}
\newcommand{\mn}{\ensuremath{\mathbb N}}
\newcommand{\mq}{\ensuremath{\mathbb Q}}
\newcommand{\half}{\tfrac{1}{2}}
\newcommand{\f}{f\times \chi}
\newcommand{\summ}{\mathop{{\sum}^{\star}}}
\newcommand{\chiq}{\chi \bmod q}
\newcommand{\chidb}{\chi \bmod db}
\newcommand{\chid}{\chi \bmod d}
\newcommand{\sym}{\text{sym}^2}
\newcommand{\hhalf}{\tfrac{1}{2}}
\newcommand{\sumstar}{\sideset{}{^*}\sum}
\newcommand{\sumprime}{\sideset{}{'}\sum}
\newcommand{\sumprimeprime}{\sideset{}{''}\sum}
\newcommand{\sumflat}{\sideset{}{^\flat}\sum}
\newcommand{\shortmod}{\ensuremath{\negthickspace \negthickspace \negthickspace \pmod}}
\newcommand{\V}{V\left(\frac{nm}{q^2}\right)}
\newcommand{\sumi}{\mathop{{\sum}^{\dagger}}}
\newcommand{\mz}{\ensuremath{\mathbb Z}}
\newcommand{\leg}[2]{\left(\tfrac{#1}{#2}\right)}
\newcommand{\muK}{\mu_{\omega}}
\newcommand{\thalf}{\tfrac12}
\newcommand{\lp}{\left(}
\newcommand{\rp}{\right)}
\newcommand{\Lam}{\Lambda_{[i]}}
\newcommand{\lam}{\lambda}
\newcommand{\af}{\mathfrak{a}}
\newcommand{\sw}{S_{[i]}(X,Y;\Phi,\Psi)}
\newcommand{\lz}{\left(}
\newcommand{\pz}{\right)}
\newcommand{\bfrac}[2]{\lz\frac{#1}{#2}\pz}
\newcommand{\odd}{\mathrm{\ primary}}
\newcommand{\even}{\text{ even}}
\newcommand{\res}{\mathrm{Res}}
\newcommand{\sumn}{\sumstar_{(c,1+i)=1}  w\left( \frac {N(c)}X \right)}
\newcommand{\lab}{\left|}
\newcommand{\rab}{\right|}
\newcommand{\Go}{\Gamma_{o}}
\newcommand{\Ge}{\Gamma_{e}}
\newcommand{\M}{\widehat}
\newcommand{\ch}{\mathrm{conv}}
\newcommand{\A}{\mathfrak{a}}

\theoremstyle{plain}
\newtheorem{conj}{Conjecture}
\newtheorem{remark}[subsection]{Remark}

\makeatletter
\def\widebreve{\mathpalette\wide@breve}
\def\wide@breve#1#2{\sbox\z@{$#1#2$}%
     \mathop{\vbox{\m@th\ialign{##\crcr
\kern0.08em\brevefill#1{0.8\wd\z@}\crcr\noalign{\nointerlineskip}%
                    $\hss#1#2\hss$\crcr}}}\limits}
\def\brevefill#1#2{$\m@th\sbox\tw@{$#1($}%
  \hss\resizebox{#2}{\wd\tw@}{\rotatebox[origin=c]{90}{\upshape(}}\hss$}
\makeatletter

\title[Shifted moments of angular Hecke $L$-functions]{Shifted moments of angular Hecke $L$-functions}

\author[P. Gao]{Peng Gao}
\address{School of Mathematical Sciences, Beihang University, Beijing 100191, China}
\email{penggao@buaa.edu.cn}

\begin{abstract}
In this paper, we establish upper bounds for shifted moments of angular Hecke $L$-function  over imaginary quadratic number fields of class number one under 
the generalized Riemann hypothesis. The result is used to derive bounds for moments of angular Hecke character sums. 
\end{abstract}

\maketitle

\noindent {\bf Mathematics Subject Classification (2010)}: 11M06, 11M41  \newline

\noindent {\bf Keywords}: angular Hecke $L$-functions, shifted moments, imaginary quadratic number fields

\section{Introduction}\label{sec 1}

  It is an important subject in number theory to study moments of central values of families of $L$-functions as they have many applications. 
  A method developed by K. Soundararajan in \cite{Sound2009} with its refinement by A. J. Harper \cite{Harper}  allows 
  one to obtain sharp upper bounds for these moments under the generalized Riemann hypothesis (GRH).
  
  Instead of the central values, one may also investigate moments of families of $L$-functions at points on the critical line and we call these 
  shifted moments. In \cites{Munsch17, Szab},  shifted moments of Dirichlet $L$- functions to a fixed modulus were investigated under GRH and the results were 
  applied to achieve bounds for moments of the corresponding Dirichlet character sums. 
 
   Motivated by the results in \cites{Munsch17, Szab}, it is the aim of this paper to first study the shifted moments of angular Hecke $L$- functions and then to apply our 
   result to obtain bounds for moments of the corresponding angular character sums. 
   
   In this paper, we focus on imaginary quadratic number fields $K$ of class number one. It is well-known that
 (see  \cite[(22.77)]{iwakow}) such a $K$ is given by $K=\mq(\sqrt{d})$ with
\begin{align}
\label{dvalue}
 d \in \{-1, -2, -3, -7,-11,-19,-43,-67,-163 \}.
\end{align}
 
   We  denote $\A$ an integral ideal of the ring of integers $\mathcal O_K$ of $K$. 
  As $\mathcal O_K$ is a principal ideal domain, the generators of $\A=(\alpha ) \subseteq \mathcal O_K$ differ
  by multiplication of an element in the group of units  $U_K$ of $K$. We shall hence fix the generator $\alpha$ of $(\alpha )$  so that
  $\theta_{\A} := \theta_{\alpha}$ lies in $[-\pi /|U_K|, \pi/|U_K|)$,  where $|U_K|$ is the order of $U_K$ and $\theta_{\alpha}$ denotes the argument of $\alpha$.

  Recall (see \cite[Section 2.1]{Holm25}) that  $|U_K|$ equals $4$ when $d=-1$, equals $6$ when $d=-3$ and equals $2$ for other values of $d$ listed
 in \eqref{dvalue}.  In what follows, we denote $M_K$ the integer such that $M_K=|U_K|/2$.

  We define the angular Hecke character $\lambda_{M_Kl}$ for any $l \in \mz$ such that 
\begin{equation*}
\lambda_{M_Kl}(\A)= \left(\frac {\alpha}{ \overline{\alpha}}\right)^{M_Kl} = e^{2i M_Kl \theta_{\A}}.
\end{equation*}
  The presence of $M_K$ ensures that $\lambda_{M_Kl}$ is trivial on $U_K$, and therefore it defines a Hecke character of trivial infinite type. 
  
  Denote $N(\A)$ the norm of the integral ideal $\A$. Throughout the paper, we reserve the letter $p$ for a rational prime number and $\varpi$ 
  for a prime in $K$, by which we mean the ideal $(\varpi)$ generalized by $\varpi$ is a prime ideal. For each $\lambda_{M_Kl}$,  
  the associated Hecke $L$-function is defined to be
\begin{align}
\label{Lexp}
L(s, \lambda_{M_Kl}) &= \ \sum_{\substack{\A \subseteq \mathcal O_K\\ \A \neq 0}}\frac{\lambda_{M_Kl}(\A)}{N(\A)^{s}}=\prod_{(\varpi) }\left(1-\frac{\lambda_{M_Kl}(\varpi)}{N(\varpi)^s}\right)^{-1}, \hspace{5mm}
\Re(s)>1.
\end{align}
  
  In \cite{Gao2026-10}, the author obtained asymptotical formulas for the first and second moments of the families of angular Hecke $L$-functions at the central
   point. In this paper, we shall apply the method in \cite{Sound2009} and \cite{Harper} to 
  establish upper bounds for shifted moments of these angular Hecke $L$-functions at points on the critical line. 
  Denote $D_K$ the discriminant of $K$ and write $\chi^{(D_K)}$ for $\leg {D_K}{\cdot}_{\mz}$, where $\leg {\cdot}{\cdot}_{\mz}$ is the Kronecker symbol in $\mz$. Our result is as follows. 
\begin{theorem}
\label{t1}
 With the notation as above and assuming the truth of GRH. Let $K$ be any imaginary quadratic number field of class number one. Let $M_K$ be the integer such that $2M_K=|U_K|$. Let $k\geq 1$ be a fixed integer and $a_1,\ldots, a_{k}$, $A$ be fixed positive real numbers. Suppose that $X$ is a large real number and $t=(t_1,\ldots ,t_{k})$ a real $k$-tuple with $|t_j|\leq  X^A$. Then
\begin{align}
\label{Lprodbounds}
\begin{split}
  & \sum_{1 \leq l \leq X} \big| L\big( \tfrac12+it_1, \lambda_{M_Kl} \big ) \big|^{a_1} \cdots \big| L\big( \tfrac12+it_{k}, \lambda_{M_Kl}  \big) \big|^{a_{k}} \\
 \ll &   X(\log X)^{3(a_1^2+\cdots +a_{k}^2)/8} \\
  & \times \prod_{1\leq j<m \leq k} \Big|\zeta \Big(1+i(t_j-t_m)+\frac 1{\log X} \Big) \Big|^{3a_ja_m/4}\Big|\zeta \Big(1+i(t_j+t_m)+\frac 1{\log X} \Big) \Big|^{a_ja_m/4} \\
  & \times \prod_{1\leq j\leq k} \Big|\zeta \Big(1+2it_j+\frac 1{\log X} \Big) \Big|^{a^2_j/8+a_j/2} \\
  & \times \prod_{1\leq j<m \leq k} \Big|L \Big(1+i(t_j-t_m)+\frac 1{\log X}, \chi^{(D_K)} \Big) \Big|^{3a_ja_m/4}\Big|L\Big(1+i(t_j+t_m)+\frac 1{\log X}, \chi^{(D_K)} \Big) \Big|^{a_ja_m/4} \\
  & \times \prod_{1\leq j\leq k} \Big|L \Big(1+2it_j+\frac 1{\log X}, \chi^{(D_K)} \Big) \Big|^{a^2_j/8-a_j/2}. 
\end{split}
\end{align}
 Here the implied constant depends on $k$, $A$ and the $a_j$'s, but not on the $t_j$'s.
\end{theorem}

   Following the approach in \cite{Szab}, we may now apply Theorem \ref{t1} to obtain bounds for the moments of angular Hecke character sum given by
\begin{align}
\label{charsum}
 \sum_{1 \leq l \leq X}\Big | \sum_{\substack{ (n) \\ N(n) \leq Y}}\lambda_{M_Kl} (n)\Big |^{2m}, 
\end{align}
 where $m$, $X$, and $Y$ are positive real numbers. 
 
  For this, we need 
   to majorize the Riemann zeta functions as well as $L(s,  \chi^{(D_K)})$ appearing in \eqref{Lprodbounds}.  
   The estimations given in \eqref{mertenstype} and \eqref{mertenstypesympower} in Section \ref{sec2.1} below allow us to deduce the following consequence of Theorem \ref{t1}.
\begin{corollary}
\label{cor1quad}
With the notation as above and assumptions as Theorem~\ref{t1}, we have for $a_j \geq 4, 1 \leq j \leq k$, 
\begin{align*}
\begin{split}
 & \sum_{1 \leq l \leq X} \big| L\big( \tfrac12+it_1, \lambda_{M_Kl} \big ) \big|^{a_1} \cdots \big| L\big( \tfrac12+it_{k}, \lambda_{M_Kl}  \big) \big|^{a_{k}} \\
 \ll & X(\log X)^{3(a_1^2+\cdots +a_{k}^2)/8} \prod_{1\leq j<m\leq k} g_1(|t_j-t_m|)^{3a_ja_m/4}g_1(|t_j+t_m|)^{a_ja_m/4}\prod_{1\leq i\leq k} g_1(|2t_j|)^{a^2_j/8+a_j/2} \\
 & \times \prod_{1\leq j<m\leq k} g_2(|t_j-t_m|)^{3a_ja_m/4}g_2(|t_j+t_m|)^{a_ja_m/4}\prod_{1\leq j\leq k} g_2(|2t_j|)^{a^2_j/8-a_j/2}, 
\end{split}
\end{align*}
  where $g_1:\mathbb{R}_{\geq 0} \rightarrow \mathbb{R}$ is defined by
\begin{equation}
\label{g3Def}
g_1(x) =\begin{cases}
\log X,  & \text{if } x\leq 1/\log X  \text{ or } x > e^X, \\
1/x, & \text{if }   1/\log X < x\leq 10, \\
\log \log x, & \text{if }  10 < x \leq e^{X},
\end{cases}
\end{equation}
  and where $g_2:\mathbb{R}_{\geq 0} \rightarrow \mathbb{R}$ is defined by
\begin{equation}
\label{g4Def}
g_2(x) =\begin{cases}
1,  & \text{if } x \leq e^e, \\
\log \log x, & \text{if }   e^e < x \leq e^{X}, \\
 \log X , & \text{if }  x > e^{X}.
\end{cases}
\end{equation}
 Here the implied constant depends on $k$, $A$ and the $a_j$'s, but not on $X$ or the $t_j$'s.
\end{corollary} 
 
    We consider instead in this paper a smoothed version of the sum given \eqref{charsum}. Let $\Phi(t)$ be a fixed non-negative smooth function supported on $(0,1)$ with $\Phi(t)=1$ for $t \in (1/X^{\varepsilon_0}, 1-1/X^{\varepsilon_0})$ and $\Phi^{(j)}(t) \ll_j X^{j\varepsilon_0}$ for all integers $j \geq 0$.  Here $0<\varepsilon_0<1$ is a fixed constant. We define
\begin{align*}
 S_{m}(X,Y;K, \Phi) :=  \sum_{1 \leq l \leq X}\Big | \sum_{ \substack{ (n)}}\lambda_{M_Kl} (n)\Phi\Big(\frac {N(n)}{Y}\Big)\Big |^{2m}.
\end{align*}    
   
    With the help of Corollaries~\ref{cor1quad}, we shall derive the following bound for $S_{m}(X,Y;K, \Phi)$ in Section \ref{sec: mainthmquad}.
\begin{theorem}
\label{quadraticmean}
With the notation as above and assuming the truth of GRH, we have, for large $X$, $Y$ and $m \geq 4$,
\begin{align}
\label{mainestimationmlarge}
 S_{m}(X,Y;K, \Phi)  \ll XY^m(\log X)^{2m^2-m+4}.
\end{align}
\end{theorem}

\section{Preliminaries}
\label{sec 2}

\subsection{Splitting of rational primes in imaginary quadratic number fields}
\label{sect: Kronecker}

   Let $K=\mq(\sqrt{d})$ be an imaginary quadratic number field of class number one, where $d$ takes the value given in \eqref{dvalue}. 
   The following facts can be found in \cite[Section 3.8]{iwakow}. Recall that $D_K$ denotes the discriminant of $K$, then we have
\begin{align*}
   D_K & =\begin{cases}
     d \qquad & \text{if $d \equiv 1 \pmod 4$}, \\
     4d \qquad & \text{if $d \equiv 2, 3 \pmod 4$}.
    \end{cases}
\end{align*}

   Recall the Dirichlet character $\chi^{(D_K)}$ from Section \ref{sec 1}. Let $\mathfrak{p}$ denote a prime ideal in $\mathcal O)K$.  A rational prime ideal $(p) \in \mz$ in $\mathcal O_K$ factors in the following way:
\begin{align}
\label{splitcriteria}
\begin{split}
    & \chi^{(D_K)}(p) =0, \text{then $p$ ramifies}, \ (p)=\mathfrak{p}^2 \quad \text{with} \quad N(\mathfrak{p})=p, \\
    &  \chi^{(D_K)}(p)  =1, \text{then $p$ splits}, \ (p)=\mathfrak{p}\overline{\mathfrak{p}} \quad \text{with} \quad  \mathfrak{p} \neq \overline{\mathfrak{p}}
    \quad \text{and} \quad  N(\mathfrak{p})=p, \\
    &  \chi^{(D_K)}(p) =-1, \text{then $p$ is inert}, \ (p)=\mathfrak{p} \quad \text{with} \quad  N(\mathfrak{p})=p^2.\
\end{split}
\end{align}

\subsection{Sums over primes}
\label{sec2.1}

  We note the following result that evaluates various sums over primes asymptotically .
\begin{lemma}
\label{RS}
 With the notation as above and let $x \geq 2, \alpha \geq 0$. Let $K$ be an algebraic number field and let $\chi^{(D_K)}$ be the Kronecker symbol $\leg {D_K}{\cdot}_{\mz}$. We have, for some constant $B_K$,
\begin{align}
\label{mertenpartialsummation}
\sum_{\substack{ (\varpi) \\ N(\varpi) \le x}}  \frac {\log N(\varpi)}{N(\varpi)} =& \log x + O(1),  \\
\label{merten}
\sum_{\substack{ (\varpi) \\ N(\varpi) \le x}} \frac{1}{N(\varpi)} =& \log \log x + B_K+ O_K\Big(\frac{1}{\log x}\Big), \\
\label{mertenstype}
  \sum_{p\leq x} \frac{\cos(\alpha \log p) }{p}=& \log |\zeta(1+1/\log x+i\alpha)| +O(1)
  \leq
\begin{cases}
\log\log x+O(1)            & \text{if }  \alpha \leq 1/\log x \text{ or } \alpha> e^x  ,   \\
\log(1/\alpha)+O(1)        & \text{if }  1/\log x< \alpha \leq 10,   \\
\log\log\log \alpha + O(1) & \text{if }   10 < \alpha \leq e^x.
\end{cases} \\
\label{mertenstypesympower}
  \sum_{p\leq x} \frac{\cos(\alpha \log p)\chi^{(D_K)}(p) }{p}=& \log |L(1+1/\log x+i\alpha, \chi^{(D_K)})| +O(1)
  \leq 
\begin{cases}
 O(1)            & \text{if }  \alpha \leq e^e,   \\
\log\log\log \alpha + O(1)        & \text{if }  e^e < \alpha \leq e^{x},   \\
 \log\log x & \text{if }  \alpha > e^{x}.
\end{cases} 
\end{align}
The estimates in the first two cases of \eqref{mertenstype} are unconditional and the third holds under the Riemann hypothesis. The estimates in the first and last case of \eqref{mertenstypesympower} are unconditional and the second one holds under GRH.

 Moreover, denote $\chi_0$ the primitive principal Hecke character and let $\chi$ be a primitive Hecke character modulo $\mathfrak a$. 
 Assume that GRH holds for $L(s, \chi)$, then we have
\begin{align}
\label{PIT}
 \sum_{\substack{ (\varpi) \\ N(\varpi) \le x}} \chi(\varpi) \log N(\varpi) = & \delta_{\chi=\chi_0}x+O(\sqrt{x} \left(\log 2N(\mathfrak a)x)^2 \right),
\end{align}
 where $\delta_{\chi=\chi_0}=1$ if $\chi=\chi_0$ and $\delta_{\chi=\chi_0}=0$ otherwise.
\end{lemma}
\begin{proof}
 The estimations given in \eqref{mertenpartialsummation} and \eqref{merten} can be found in Lemma 2.3 and Lemma 2.4 in \cite{Rosen99}.  The equality in  \eqref{mertenstype} is a special case given in \cite[Lemma 3.2]{Kou} and the other estimations can be found in \cite[Lemma 2]{Szab}. The proof of \eqref{mertenstypesympower} follows by arguing similar to those given in the proof of \cite[Lemma 2.6]{G&Zhao24-12}, upon using \cite[Lemma 2.5]{G&Zhao24-12}.
 The estimation given in \eqref{PIT} follows from \cite[(5.55)]{iwakow}. 
\end{proof}

\subsection{Angular Hecke character sums}

 The following result provides an estimation concerning the average of angular Hecke character sums.
\begin{lemma}
\label{prsum}
 With the notation as above, we have for $\A \neq (0)$, 
\begin{align}
\label{sumlambda}
\begin{split}
 \sum_{1 \leq l\leq X}\lambda_{M_Kl}(\A) 
=&
\begin{cases}
 X +O(1),  \ \text{if $\A=(n)$ with $n \in \mn$},  \\
 O(N(\A)^{1/2}), \ \text{otherwise}. 
\end{cases}
\end{split}
\end{align}
\end{lemma}
\begin{proof}
  The case when $\A=(n)$ with $n \in \mn$ is trivial as in this case, we have $\lambda_{M_Kl}(\A)=1$ for all $l$. Otherwise, we deduce from \cite[(3.10)--(3.12)]{Gao2026-10} that the desired estimation holds. This completes the proof of the lemma. 
\end{proof}

   The above lemma allows us to establish a mean value estimate similar to \cite[Lemma 3]{Sound2009}. 
\begin{lemma}
\label{lem:2.5}
With the notation as above, let $X$ and $y$ be real numbers and $m$ a positive integer.  For fixed $\varepsilon$ with $0<\varepsilon<1$ and any complex numbers $a(\varpi)$, we have 
 \begin{align*}
   &  \sum_{\substack{1\leq l \leq X}}\left|\sum_{\substack{ (\varpi) \\ N(\varpi) \leq y}}\frac{a(\varpi)\lambda_{M_kl}(\varpi)}{N(\varpi)^{1/2}}\right|^{2m} \\
  \ll_{\varepsilon} & X\sum^{m}_{i=0} \binom {m}{i}^2 i! \left(\sum_{\substack{(\varpi) \\ N(\varpi) \leq y}} \frac {|a(\varpi)|^2}{N(\varpi)}\right)^{2i} \left(\sum_{p \leq y} \frac {|a(p)|}{p}\right)^{2m-2i}+y^{2m} \left( \sum_{\substack{(\varpi) \\ N(\varpi)\leq y}} \frac {|a(\varpi)|^2}{N(\varpi)} \right)^{m}.
 \end{align*}
\end{lemma}
\begin{proof}
  Note that
\begin{align*}
  \sum_{\substack{1\leq l \leq X}}  \left|\sum_{\substack{ (\varpi) \\ N(\varpi) \leq y}}\frac{a(\varpi)\lambda_{M_kl}(\varpi)}{N(\varpi)^{1/2}}\right|^{2m} = 
& \sum_{\substack{1\leq l \leq X}}  \left| \sum_{\substack{(\varpi_i) \\ 1 \leq i \leq m \\ N(\varpi_i)\leq y}} \frac{a(\varpi_1) \dots a(\varpi_m) \ }{\sqrt{N(\varpi_1) \dots N(\varpi_m)}}
\lambda_{M_kl}(\varpi_1\cdots \varpi_m) \right|^2. 
\end{align*}
  Expanding out the square in the last sum above and then applying Lemma \ref{prsum} to evaluate the resulting sums, we see that the last expression above is
\begin{align}
\label{2kbound}
\ll &
X \sum_{\substack{(\varpi_i) \\ 1 \leq i \leq 2m \\  N(\varpi_i)\leq y \\ \varpi_1 \dots \varpi_m \dots \overline{\varpi_{m+1}} \dots \overline{\varpi_{2m}} \in \mn }} \frac{|a(\varpi_1) \dots a(\varpi_{2m}) |}{\sqrt{N(\varpi_1 \dots \varpi_{2m})}} + O\Big(\sum_{\substack{(\varpi_i) \\ 1 \leq i \leq 2m \\ N(\varpi_i)\leq y}}|a(\varpi_1) \dots a(\varpi_{2m})|   \Big).
\end{align}
 To estimate the first term above, we observe that $\varpi_1 \dots \varpi_m \dots \overline{\varpi_{m+1}} \dots \overline{\varpi_{2m}} \in \mn$  precisely when there is a way to pair the primes in $\{ \varpi_1, \cdots, \varpi_m \}$ that are not in $\mn$ with those in $\{ \varpi_{m+1}, \cdots, \varpi_{2m}\}$. Suppose there are $i$ such primes in  $\{ \varpi_1, \cdots, \varpi_m \}$ that are not in $\mn$ with $1 \leq i \leq m$, then the number of such pairing equals 
\begin{align*}
  \binom {m}{i}^2 i!.
\end{align*}   
   Note further that when $\varpi=\varpi'$, we have $a(\varpi)a(\varpi')/\sqrt{N(\varpi\varpi')} \ll |a(\varpi)|^2/N(\varpi)$. Moreover, if $\varpi=p$ is a rational prime, then we have $a(\varpi)/\sqrt{N(\varpi)}=a(p)/p$. It follows that we have
\begin{align}
\label{maintermcharbound}
 & \sum_{\substack{(\varpi_i) \\ 1 \leq i \leq 2m \\  N(\varpi_i)\leq y \\ \varpi_1 \dots \varpi_m \dots \overline{\varpi_{m+1}} \dots \overline{\varpi_{2m}} \in \mn }} \frac{|a(\varpi_1) \dots a(\varpi_{2m}) |}{\sqrt{N(\varpi_1 \dots \varpi_{2m})}} \ll   \sum^{m}_{i=0} \binom {m}{i}^2 i! \left(\sum_{\substack{(\varpi) \\ N(\varpi) \leq y}} \frac {|a(\varpi)|^2}{N(\varpi)}\right)^{2i} \left(\sum_{p \leq y} \frac {|a(p)|}{p}\right)^{2m-2i}.
\end{align}

   On the other hand, the Cauchy-Schwarz inequality yields
\begin{align}
\label{errortermcharbound}
\begin{split}
\sum_{\substack{(\varpi_i) \\ 1 \leq i \leq 2m \\ N(\varpi_i)\leq y}}|a(\varpi_1) \dots a(\varpi_{2m})|  \ll &  \left( \sum_{\substack{(\varpi) \\ N(\varpi)\leq y}}|a(\varpi)| \right)^{2m} \ll  \left( \sum_{\substack{(\varpi) \\ N(\varpi)\leq y}} \frac {|a(\varpi)|^2}{N(\varpi)} \right)^{m}\left(\sum_{\substack{(\varpi) \\ N(\varpi)\leq y}} N(\varpi) \right)^{m}  \\
\ll  y^{2m} \left( \sum_{\substack{(\varpi) \\ N(\varpi)\leq y}} \frac {|a(\varpi)|^2}{N(\varpi)} \right)^{m}.
\end{split}
\end{align}
 We readily deduce the assertion of the lemma  from \eqref{2kbound}--\eqref{errortermcharbound}. This completes the proof of the lemma.
\end{proof}

\subsection{Upper bound for $\log |L(\frac 12+it, \lambda_{M_Kl})|$ }
 
  Our next result provides an upper bound of $\log |L(1/2+it, \Lambda_{M_Kl})|$ in terms of a sum involving prime powers.
\begin{prop}
\label{lem: logLbound}
With the notation as above and let $l \geq 1$. Assume the truth of GRH for $L(s, \lambda_{M_Kl})$. Let $x \geq 2$ and let $\eta_0=0.4912\ldots$
denote the unique positive real number satisfying $e^{-\eta_0} = \eta_0+\frac {\eta^2_0}{2}$.  Then, for any $\eta \geq \eta_0$,
\begin{align}
\label{logLupperbound}
\begin{split}
\log |L(\frac 12+it, \lambda_{M_Kl})| \le &  \Re  \sum_{\substack {(\varpi^m) \\ N(\varpi^m) \leq x \\  m \geq 1}} \frac {\log N(\varpi)  \cdot  \lambda_{M_Kl}(\varpi^{m})}{mN(\varpi)^{m(\sigma_0+it)}}\cdot
\frac {\log \big (\frac {x}{N(\varpi^m)}\big )}{\log x} \\
& + (1+\eta)(\log (l+2)+\log (|t|+2))+ O\Big( \frac{\eta +1}{\log x}+1\Big).
\end{split}
\end{align}
\end{prop}
\begin{proof}
We write $s=\sigma+it$ and interpret $\log |L(s, \lambda_{M_Kl})| $ as $-\infty$ when $L(s, \lambda_{M_Kl})=0$. We may
therefore assume that $L(s, \lambda_{M_Kl}) \neq 0$ in the rest of the proof. It is shown in \cite[Section 2]{Waxman21} that for $l \neq 0$, the complete Hecke $L$-function
\begin{equation}
\label{fcneqn}
\Lambda_{l}(s) = |D_K|^{s/2}(2\pi)^{-s}\Gamma(s+|M_Kl|)L(s, \lambda_{M_Kl})
\end{equation}
  is entire and have an analytic continuation to the whole complex plane that satisfies the functional equation
\begin{equation*}
 \Lambda_{l}(s) = \Lambda_{l}(1-s).
\end{equation*} 

  As $\Gamma(s)$ has simple poles at the non-positive rational integers
(see \cite[\S 10]{Da}), the above expression of $\Lambda_l(s)$ implies that $L(s,\lambda_{M_Kl})$ has simple zeros at $s=-k-M_Kl$ for all $k \geq 0$, 
and we call these the trivial zeros of $L(s,\lambda_{M_Kl})$.  The non-trivial zeros of $L(s,\lambda_{M_Kl})$ are precisely the zeros of $\Lambda_l(s)$. 

  Let $\rho=\tfrac 12+i\gamma$ run over the non-trivial zeros of $L(s,\lambda_{M_Kl})$ so that we have $\gamma \in \rear$ by GRH.  We then deduce from \cite[Theorem 5.6]{iwakow}
and the observation that $L(s,\lambda_{M_Kl})$ is analytic at $s=1$ when $l \geq 1$ that for constants $B_0$, $B_1=B_1(l)$,
\begin{align}
\label{Lproductzeros}
  \Lambda_l(s)=e^{B_0 + B_1s}\prod_{\rho}\left(1 - \frac{s}{\rho}\right)e^{\frac{s}{\rho}}. 
\end{align}
 
   We substitute the expression given in \eqref{fcneqn} for $\Lambda_l(s)$ in \eqref{Lproductzeros}, 
   and then take the logarithmic derivative on both sides of the resulting expression. By further taking the real parts, we obtain that
\begin{align}
\label{Lderivreal}
-\Re \frac{L'}{L}(s, \lambda_{M_Kl})=\log \left(\frac{\sqrt{|D_K|}}{2\pi} \right)+\Re \frac{\Gamma'}{\Gamma} \left(s +M_Kl\right)
- \Re B_1 - \sum_{\rho}\Re \left(\frac{1}{s-\rho} + \frac{1}{\rho}\right).
\end{align}

    On the other hand, we note that by \cite[(5.29)]{iwakow},
\begin{equation*}
\Re(B_1)=-\sum_\rho\Re(\rho^{-1}).
\end{equation*}
  It follows from the above and \eqref{Lderivreal} that we have
\begin{align}
\label{LprimeLbound}
-\Re{\frac{L'}{L}(s, \lambda_{M_Kl})} = \Re \frac{\Gamma'}{\Gamma} \left( s +M_Kl\right)-F(s) +O(1),
\end{align}
  where
\begin{align*}
 F(s) := \Re{\sum_{\rho} \frac{1}{s-\rho}}= \sum_{\rho} \frac{\sigma-1/2}{(\sigma-1/2)^2+
 (t-\gamma)^2}.
\end{align*}

  Note that $\frac {\Gamma'}{\Gamma} (s)=\log s+O(|s|^{-1})$ (see \cite[(C. 17)]{MVa1}). It follows from this and \eqref{LprimeLbound} that, for $s=\sigma+it$ with $\sigma$, $t \in \mr$ and $|\sigma|$ is bounded, 
\begin{align}
\label{LprimeLbound1}
\begin{split}
-\Re{\frac{L'}{L}(\sigma+it, \lambda_{M_Kl})} \leq & \log (l+2)+\log (|t|+2) - F(\sigma+it)+O(1).
\end{split}
\end{align}

  Upon integrating the last expression above from $\sigma=1/2$ to
$\sigma=\sigma_0 > \frac 12$, we obtain that
\begin{align}
\label{Lprimediffbound}
\begin{split}
  \log |L(\tfrac 12+it,  \lambda_{M_Kl})| - \log & |L(\sigma_0+it,  \lambda_{M_Kl})| \leq (\sigma_0-\half) \Big(\log (l+2)+\log (|t|+2)+O(1)\Big)  -\int\limits_{\half}^{\sigma_0} F(\sigma+it) \dif \sigma \\
 =&  (\sigma_0-\half) \Big (\log (l+2)+\log (|t|+2)  +O(1) \Big)-\half \sum_{\rho} \log \frac {(\sigma_0-\half)^2+(t-\gamma)^2}{(t-\gamma)^2} \\
 \leq &  (\sigma_0-\half) \Big (\log (l+2)+\log (|t|+2)-\half F(\sigma_0+it) +O(1)  \Big),
\end{split}
\end{align}
  where the inequality above follows from the observation that $\log (1+x^2) \geq x^2/(1+x^2)$. 

  Observe from \eqref{Lexp} that we have for $\Re(s)>1$, 
  \begin{align*}
\begin{split}
 -\frac{L^{\prime}}{L}(s, \lambda_{M_Kl})
 =& \sum_{\substack { (\varpi^m) \\ m \geq 1}} \frac {\log N(\varpi) \cdot \lambda_{M_Kl}(\varpi^{m})}{N(\varpi)^{ms}}.
\end{split}
\end{align*}

 We integrate both sides above to see that
$$
 \frac{1}{2\pi i} \int\limits_{(c)} -\frac{L^{\prime}}{L}(s+w, \lambda_{M_Kl})
 \frac{x^w}{w^2} \dif w  = \sum_{\substack {(\varpi^m) \\ N(\varpi^m) \leq x \\  m \geq 1}} \frac {\log N(\varpi)  \cdot  \lambda_{M_Kl}(\varpi^{l})}{N(\varpi)^{ms}}
\log \Big (\frac {x}{N(\varpi^m)}\Big ),
 $$
   where $c>2$ is a large real number. Upon shifting the line of integration in the above expression to the left and collecting residues, we see that
\begin{align*}
\frac{1}{2\pi i} \int\limits_{(c)} -\frac{L^{\prime}}{L}(s+w, \lambda_{M_Kl})
 \frac{x^w}{w^2} \dif w 
 = -\frac{L^{\prime}}{L}(s, \lambda_{M_Kl}) \log x - \Big(\frac{L^{\prime}}{L}(s, \lambda_{M_Kl})\Big)^{\prime}
 -\sum_{\rho} \frac{x^{\rho-s}}{(\rho-s)^2} -\sum_{k=0}^{\infty}\frac{x^{-k-M_Kl-s}}{(k+M_Kl+s)^2}.
 \end{align*}

  We compare the above two expressions to derive unconditionally that for any $x \ge 2$,
 \begin{align}
\label{Lprimeseries}
\begin{split}
 -\frac{L^{\prime}}{L}(s, \lambda_{M_Kl})=& \sum_{\substack {(\varpi^m) \\ N(\varpi^m) \leq x \\  m \geq 1}} \frac {\log N(\varpi)  \cdot  \lambda_{M_Kl}(\varpi^{m})}{N(\varpi)^{ms}} \cdot
\frac {\log \big (\frac {x}{N(\varpi^m)}\big )}{\log x}+ \frac{1}{\log x} \Big(\frac{L^{\prime}}{L}(s, \lambda_{M_Kl})\Big)^{\prime}
 + \frac{1}{\log x} \sum_{\rho} \frac{x^{\rho-s}}{(\rho-s)^2} \\
& + \frac{1}{\log x} \sum_{k=0}^{\infty}\frac{x^{-k-M_Kl-s}}{(k+M_Kl+s)^2}.
\end{split}
 \end{align}

 Integrating the real parts on both sides of \eqref{Lprimeseries} from $\sigma_0>1/2$ to $\infty$, we see that for $x
\geq 2$,
 \begin{align}
\label{logL}
\begin{split}
\log |L(\sigma_0+it, \lambda_{M_Kl})| =& \Re \Big( \sum_{\substack {(\varpi^m) \\ N(\varpi^m) \leq x \\  m \geq 1}} \frac {\log N(\varpi)  \cdot  \lambda_{M_Kl}(\varpi^{m})}{lN(\varpi)^{m(\sigma_0+it)}}\cdot
\frac {\log \big (\frac {x}{N(\varpi^m)}\big )}{\log x}
  - \frac{1}{\log x} \frac{L^{\prime}}{L}(\sigma_0+it, \lambda_{M_Kl})\\
  & + \frac{1}{\log x} \sum_{\rho}
 \int\limits_{\sigma_0}^{\infty} \frac{x^{\rho-s}}{(\rho-s)^2} \dif \sigma +O\Big(\frac{1}{\log x}\Big)\Big).
\end{split}
\end{align}
 Note that
$$
\sum_{\rho}\Big|\int\limits_{\sigma_0}^{\infty} \frac{x^{\rho -s}}{(\rho -s)^2} \dif \sigma\Big|
\le \sum_{\rho}\int\limits_{\sigma_0}^{\infty}\frac{ x^{1/2-\sigma}}{|\sigma_0-\rho|^2} \dif \sigma
= \sum_{\rho}\frac{x^{1/2-\sigma_0}}{|\sigma_0+it-\rho|^2 \log x}= \frac{x^{1/2-\sigma_0}F(\sigma_0+it)}{(\sigma_0-1/2)\log x}.
$$
 Applying the above and \eqref{LprimeLbound1} in \eqref{logL} leads to
 \begin{align}
\label{logLbound}
\begin{split}
 \log |L(\sigma_0+it, \lambda_{M_Kl})|
 &\le  \Re \Big(  \sum_{\substack {(\varpi^m) \\ N(\varpi^m) \leq x \\  m \geq 1}} \frac {\log N(\varpi)  \cdot  \lambda_{M_Kl}(\varpi^{m})}{mN(\varpi)^{m(\sigma_0+it)}}\cdot
\frac {\log \big (\frac {x}{N(\varpi^m)}\big )}{\log x}
 +  \frac {\log (l+2)+\log (|t|+2)}{ \log x}  \\
& \hskip .5 in+F(\sigma_0+it) \Big( \frac{x^{1/2-\sigma_0}}{(\sigma_0-1/2) \log^2 x} -\frac{1}{\log x}
 \Big)
 + O\Big(\frac{1}{\log x}+1\Big)\Big).
\end{split}
 \end{align}

 We add \eqref{Lprimediffbound} and \eqref{logLbound} to deduce that
\begin{align*}
  \log |L(\tfrac 12+it, \lambda_{M_Kl})|
 \le &  \Re  \sum_{\substack {(\varpi^m) \\ N(\varpi^m) \leq x \\  m \geq 1}} \frac {\log N(\varpi)  \cdot  \lambda_{M_Kl}(\varpi^{m})}{mN(\varpi)^{m(\sigma_0+it)}}\cdot
\frac {\log \big (\frac {x}{N(\varpi^m)}\big )}{\log x}
 +  (\log (l+2)+\log (|t|+2) ) \Big(\sigma_0 -\frac{1}{2} + \frac{1}{\log x} \Big )
 \\
 &\hskip .5 in +F(\sigma_0+it) \Big( \frac{x^{1/2-\sigma_0}}{(\sigma_0-1/2) \log^2 x} -\frac{1}{\log x}-\frac{1}{2} \Big( \sigma_0-\frac{1}{2} \Big)
 \Big) + O\Big(\frac{1}{\log x} +\sigma_0-\frac 12 \Big).
 \end{align*}
  We now set $\sigma_0 = 1/2+ \eta/\log x$ for $\eta \geq \eta_0$ to see that we may discard the term involving $F(\sigma_0+it)$ in the above inequality as its coefficient is negative. This now leads to the estimation given in \eqref{logLupperbound} and therefore completes the proof of the proposition. 
 \end{proof}

  Our next result provides a simplification of the above proposition.
\begin{corollary}
\label{lem: logLboundsimplified}
 With the notation as above, let $X>0$ be large and let $1 \leq l \leq X$. Assume the truth of GRH for $L(s, \lambda_{M_Kl})$ for any $1 \leq l \leq X$. We have for $t \leq X^A$ with $A>0$ being a fixed positive constant and any $2 \leq x \leq X$,
\begin{align}
\label{logLupperboundnonquad}
\begin{split}
\log |L(\tfrac12+it, \lambda_{M_Kl})| \leq &  \Re\sum_{\substack{ (\varpi) \\ N(\varpi) \leq x }} \frac{\lambda_{M_Kl}(\varpi)}{N(\varpi)^{1/2+it+1/\log x}}
 \frac {\log \big (\frac {x}{N(\varpi)}\big )}{\log x}+\frac 12\Re \sum_{\substack{   (\varpi) \\  N(\varpi) \leq \min (\log X, x^{1/2})   }} \frac{\lambda_{M_Kl} (\varpi^2)}{N(\varpi)^{1+2it}}
 +\frac{2(A+1)\log X}{\log x} \\
 &+ O\Big( \frac {\log X}{x^{1/2}\log x}+1 \Big).
\end{split}
\end{align}
\end{corollary}
\begin{proof}
  We set $\eta=1$ in \eqref{logLupperbound} and observe that the contribution of the terms with $m \geq 3$ is $O(1)$. It follows that
\begin{align}
\begin{split}
\label{logLupperbound1}
 \log  |L(\tfrac12+it, \lambda_{M_Kl})| 
  \leq & \Re \sum_{\substack{ (\varpi) \\ N(\varpi) \leq x }} \frac{\lambda_{M_Kl}(\varpi)}{N(\varpi)^{1/2+it+1/\log x}}
 \frac {\log \big (\frac {x}{N(\varpi)}\big )}{\log x} \\
 & +  \frac 12
 \Re\sum_{\substack{ (\varpi) \\ N(\varpi) \leq x^{1/2}   }} \frac{\lambda_{M_Kl} (\varpi^2)}{N(\varpi)^{1+2it+2/\log x}}  \frac{\log (x/N(\varpi^2))}{\log x}
 +\frac{2(A+1)\log X}{\log x} + O(1).
\end{split}
\end{align}

   Note that
\begin{align}
\begin{split}
\label{sumoverpeval}
 \sum_{\substack{ (\varpi) \\ N(\varpi) \leq x^{1/2}   }} \frac{\lambda_{M_Kl} (\varpi^2)}{N(\varpi)^{1+2it+2/\log x}}  \frac{\log (x/N(\varpi^2))}{\log x}
 =\sum_{\substack{ (\varpi) \\ N(\varpi) \leq x^{1/2}   }} \frac{\lambda_{M_Kl} (\varpi^2)}{N(\varpi)^{1+2it+2/\log x}}  
 - 2\sum_{\substack{ (\varpi) \\ N(\varpi) \leq x^{1/2}   }} \frac{\lambda_{M_Kl} (\varpi^2)}{N(\varpi)^{1+2it+2/\log x}}  \frac{\log N(\varpi)}{\log x}.
\end{split}
\end{align}

  By \eqref{mertenpartialsummation}, we have 
\begin{align}
\label{errorsumlogp}
  \sum_{\substack{ (\varpi) \\ N(\varpi) \leq x^{1/2}  }}  \frac{\lambda_{M_Kl} (\varpi^2)}{N(\varpi)^{1+2it+2/\log x}}  \frac{\log N(\varpi)}{\log x} \ll  \frac 1{\log x} \sum_{\substack{ (\varpi) \\ N(\varpi) \leq x^{1/2}  }}  \frac{\log N(\varpi)}{N(\varpi)} \ll 1.
\end{align}

Similarly,
\begin{align}
\label{sumpsquare}
\begin{split}
 \sum_{\substack{ (\varpi) \\ N(\varpi) \leq x^{1/2}  }} \frac{\lambda_{M_Kl} (\varpi^2)}{N(\varpi)^{1+2it+2/\log x}} 
 = & \sum_{\substack{ (\varpi) \\ N(\varpi) \leq x^{1/2}   }} \frac{\lambda_{M_Kl} (\varpi^2)}{N(\varpi)^{1+2it}}  +\sum_{\substack{ (\varpi) \\ N(\varpi) \leq x^{1/2} }} \lambda_{M_Kl} (\varpi^2)\Big (\frac{1}{N(\varpi^{1+2it+2/\log x})}-\frac 1{N(\varpi^{1+2it})} \Big ) \\
 =&\sum_{\substack{ (\varpi) \\ N(\varpi) \leq x^{1/2}   }} \frac{\lambda_{M_Kl} (\varpi^2)}{N(\varpi)^{1+2it}}  +O\Big ( \sum_{\substack{ (\varpi) \\ N(\varpi) \leq x^{1/2}   }} \frac{\log N(\varpi)}{N(\varpi) \log x}\Big ) 
 = \sum_{\substack{ (\varpi) \\ N(\varpi) \leq x^{1/2}   }} \frac{\lambda_{M_Kl} (\varpi^2)}{N(\varpi)^{1+2it}} +O(1). 
\end{split}
\end{align}

 We apply \eqref{sumoverpeval}--\eqref{sumpsquare} in \eqref{logLupperbound1} to see that
\begin{align}
\label{logLupperboundnonquad1}
\begin{split}
\log |L(\tfrac12+it, \lambda_{M_Kl})| \leq &  \Re\sum_{\substack{ (\varpi) \\ N(\varpi) \leq x }} \frac{\lambda_{M_Kl}(\varpi)}{N(\varpi)^{1/2+it+1/\log x}}
 \frac {\log \big (\frac {x}{N(\varpi)}\big )}{\log x}+\frac 12\Re \sum_{\substack{   (\varpi) \\  N(\varpi) \leq x^{1/2}  }} \frac{\lambda_{M_Kl} (\varpi^2)}{N(\varpi)^{1+2it}}
 +\frac{2(A+1)\log X}{\log x} \\
 &+ O\Big( \frac {\log X}{x^{1/2}\log x}+1 \Big).
\end{split}
\end{align}

 Now, to establish \eqref{logLupperboundnonquad}, we may further assume that $\log X \leq x^{1/2}$ and we observe that the character $\lambda_{M_Kl}^2$ 
 is also a primitive Hecke character modulo the ideal $(1)$. We then apply \eqref{PIT} and partial summation to see that
\begin{align}
\label{sumplarge}
 \sum_{\substack{(\varpi) \\ (\log X)^6 < N(\varpi) \leq x^{1/2} }}  \frac{\lambda_{M_Kl}(\varpi^2)}{N(\varpi)^{1+2it}}, \sum_{\substack{(\varpi) \\ (\log X)^6 < N(\varpi) \leq x^{1/2} }}  \frac{\lambda_{M_Kl}(\varpi^2)}{N(\varpi)^{1+2it}}   \ll 1.
\end{align}
 Moreover, \eqref{merten} leads to
\begin{align}
\label{sumpsmall}
  \sum_{\substack{(\varpi) \\ \log X \leq N(\varpi)  \leq (\log X)^6}}  \frac{\lambda_{M_Kl}(\varpi^2)}{N(\varpi)^{1+2it}} \ll \sum_{\substack{(\varpi) \\ \log X \leq N(\varpi)  \leq (\log X)^6}}  \frac{1}{N(\varpi)}  \ll 1.
\end{align}

 We apply \eqref{sumoverpeval}--\eqref{sumpsmall} in \eqref{logLupperboundnonquad1} to see that the estimation given in \eqref{logLupperboundnonquad} follows. This completes the proof of the corollary.
\end{proof}

 As a straightforward consequence of Corollary \ref{lem: logLboundsimplified}, the next result provides an upper bound for sums $\log |L(\frac 12+it,\lambda_{M_Kl} )|$ over various $t$.
\begin{corollary}
\label{lemmasum}
  With the notation as above, let $X>0$ be large and let $1 \leq l \leq X$. Assume the truth of GRH for $L(s, \lambda_{M_Kl})$ for any $1 \leq l \leq X$. Let $k$ be a positive integer and $A,a_1,a_2,\ldots, a_{k}$ fixed positive real constants.  Set $a:=a_1+\cdots+ a_{k}$.  Let $t_1,\ldots, t_{k}$ be fixed real numbers with $|t_i|\leq X^A$. For any $c \in \mathcal O_K$, let
\begin{align}
\label{hdef}
\begin{split}
  h(c) :=\frac{1}{2}  \sum^{k}_{m=1}a_mN(c)^{-it_m} .
\end{split}
\end{align}
Then for $x \geq 2$,
\begin{align*}
\begin{split}
 \sum^{k}_{m=1}  a_m  \log |L(\tfrac12+it_m, \lambda_{M_Kl} )|  \leq &  2 \Re\sum_{\substack{ (\varpi) \\ N(\varpi) \leq x }} \frac{h(\varpi)\lambda_{M_Kl}(\varpi)}{N(\varpi)^{1/2+1/\log x}}
 \frac {\log \big (\frac {x}{N(\varpi)}\big )}{\log x} \\
 & +\Re \sum_{\substack{   (\varpi) \\  N(\varpi) \leq \min (\log X, x^{1/2})   }} \frac{h(\varpi^2)\lambda_{M_Kl} (\varpi^2)}{N(\varpi)}+2(A+1)a\frac{\log X}{\log x}+ O\Big(\frac {\log X}{x^{1/2}\log x}+1\Big).
\end{split}
\end{align*}
\end{corollary}

\subsection{A first bound on shifted moments of angular $L$-functions}

   In this section, we aim to establish a crude upper bound on shifted moments of angular $L$-functions. As a preparation, we denote $\mathcal{N}_t(V, X)$ the number of angular Hecke characters $\lambda_{M_Kl}$ with $1 \leq l \leq X$ and $\log|L(\frac 12+it, \lambda_{M_Kl})|\geq V$ with $|t| \leq X^A$ for some constant $A>0$.
  We have the following upper bounds for $\mathcal{N}_t(V, X)$.
\begin{prop}
\label{propNbound}
 With the notation as above and assuming the truth of GRH for $L(s,\lambda_{M_Kl})$ for any $1 \leq l \leq X$. Let $k>0$ be a fixed real number.
 Then for $X$ large enough, 
 \begin{align*}
  \mathcal{N}_t(V,X)\ll_{k,A} 
 \left\{
 \begin{array}
  [c]{ll}
X(\log \log X)^2\exp \left(-\frac{V^2}{10^4\log \log X} \right), \quad 10\sqrt{\operatorname{log}\operatorname{log} X}\leq V \leq 10^{4+4k} \log \log X, \\ \\
 XV^2\exp\left(-(2+4k) V \right), \quad  V > 10^{4+4k} \log \log X .
\end{array}
\right .
 \end{align*}
\end{prop}
\begin{proof} 
 Let $C =\log\log\log X/10^{6+4k}$ and $z=x^{1/\log \log X}$. We set $x=X^{C/V}$ in \eqref{logLupperboundnonquad} by
    writing $M_1$ for the real part of the first sum in \eqref{logLupperboundnonquad} truncated to $N(\varpi) \leq z$ and $M_2$ 
    the complementary sum over $z < N(\varpi) \leq x$.  It  follows from this and \eqref{merten} that
\[
 \log  |L(\frac 12+it, \lambda_{M_Kl})|
    \leq M_1 + M_2 + 2(A+1)\frac {V}{C} + O(\log \log \log X).
\]
Thus if  $\log  |L(\frac 12+it, \lambda_{M_Kl})| \geq V $, then we have either
\[ M_2 \geq \frac{(A+1)V}{4C} \quad \mbox{or} \quad M_1 \geq V_1 := V \left(1-\frac{3(A+1)}{2C} \right). \]

 We now set
\[ \operatorname{meas}(X;M_1) = \# \{1 \leq l \leq X: M_1 \geq V_1 \}, \ \operatorname{meas}(X;M_2) = \# \left\{ 1 \leq l \leq X: M_2 \geq \frac{(A+1)V}{4C}  \right\}. \]

 Denote $\lceil x \rceil = \min \{ n \in \intz : n \geq x\}$ for any real $x$ and we set $m = \lceil (A+1)V/(4C) \rceil $ in Lemma \ref{lem:2.5} to deduce that
\begin{align} \label{M1bound}
\begin{split}
 \Big(\frac {(A+1)V}{4C} \Big)^{2m} \operatorname{meas}(X;M_2) \leq & \sum_{\substack{1 \leq l \leq X}} |M_2|^{2m}  \\
 \ll & X\sum^{m}_{i=0}\binom {m}{i}^2 i! \Big (\sum_{\substack{(\varpi) \\ z<N(\varpi) \leq x}} \frac {1}{N(\varpi)}\Big )^{2i}\Big (\sum_{z<p \leq x} \frac {1}{p} \Big )^{2m-2i}+x^{2m} \Big( \sum_{\substack{(\varpi) \\ z<N(\varpi) \leq x}} \frac {1}{N(\varpi)} \Big)^{m}.
\end{split}
\end{align}

    We note the following estimations concerning with factorials. 
\begin{align}
\label{Stirling}
 \binom {n}{m} \leq \left( \frac {en}{m} \right)^m, \quad \left( \frac ne \right)^n \leq n! \leq n \left( \frac ne \right)^n.
\end{align}

   We apply \eqref{Stirling} to see that for $1 \leq i \leq m$,
\begin{align*}
  \binom {m}{i}^2 i!  \leq   \Big (\frac {em}{i} \Big )^{2i} i \Big (\frac {i}{e} \Big )^{i}.
\end{align*}
  As the function $x \mapsto m^{2x}/x^{x}$ is increasing for $m \geq e$, we see that for all $i \geq 0$ and $m \geq e$,
\begin{align}
\label{binombound}
  \binom {m}{i}^2 i! \leq  e^mm^{m+1}.
\end{align}

  We deduce from this, \eqref{merten} and \eqref{M1bound} that
\begin{align*}
 \Big(\frac {(A+1)V}{4C} \Big)^{2m} \operatorname{meas}(X;M_2)
 \ll &  X m^2 (em)^{m}  \left (\left(\sum_{\substack{(\varpi) \\ z<N(\varpi) \leq x}} \frac {1}{N(\varpi)}\right)^{2m}+  \left(\sum_{\substack{z<p  \leq x}}
\frac{1}{p}\right)^{2m}\right ) \\
\ll & X m^2 (em)^{m}   \Big (\log \log \log X+O(1) \Big )^{2m}.
\end{align*}

  The above implies that
\begin{equation}
\label{equ:bd-S-2}
\operatorname{meas}(X;M_2) \ll X\operatorname{exp}\left(-\frac{V}{20C}\log V\right)\ll
\left\{
 \begin{array}
  [c]{ll}
X\exp \left(-\frac{V^2}{10^4\log \log X} \right), \quad V \leq 10^{4+4k} \log \log X, \\ \\
 X\exp\left(-(2+4k) V \right), \quad V > 10^{4+4k} \log \log X.
\end{array}
\right .
\end{equation}

  Similarly, 
  we apply Lemma \ref{lem:2.5}, \eqref{merten} and \eqref{binombound} to see that
\begin{align*}
\begin{split}
V^{2m}_1\operatorname{meas}(X;M_1) \leq  &  \sum_{1 \leq l \leq X} |M_1|^{2m} \ll X\sum^{m}_{i=0}\binom {m}{i}^2 i! \Big (\sum_{\substack{(\varpi) \\ N(\varpi)< z}} \frac {1}{N(\varpi)}\Big )^{2m-2i}\Big (\sum_{p \leq z} \frac {1}{p}\Big )^{2i}\ll  X m^2 (em)^m \Big (\log \log X \Big )^{m}.
\end{split}
\end{align*}
  We set $m=\lceil V^2_1/e^2\log \log X \rceil $ to see from the above that
\begin{align}
\label{M1bound1}
\operatorname{meas}(X ; M_1)
\ll X m^2 \left(\frac{em \log \log X}{V_1^2} \right)^m \ll  Xm^2\operatorname{exp}\left(-\frac{V_1^2}{e^3\log \log X} \right)\ll  Xm^2\operatorname{exp}\left(-\frac{V^2}{10^4\log \log X} \right).
\end{align}

  Note that we have $m \ll \log \log X$ when $V \leq 10^{4+4k}(\log \log X)$ and we have $m \ll V$ when $V > 10^{4+4k}(\log \log X)$.  
 Also, we have $\frac{V^2}{10^4\log \log X} \geq -(2+4k) V $ when $V \geq 10^{4+4k}(\log \log X)$. It follows that 
\begin{align}
\label{M1bound2}
\begin{split}
\operatorname{meas}(X ; M_1)
\ll
\left\{
 \begin{array}
  [c]{ll}
X(\log \log X)^2\exp \left(-\frac{V^2}{10^4\log \log X} \right), \quad V \leq 10^{4+4k} \log \log X, \\ \\
 XV^2\exp\left(-(2+4k) V \right), \quad V > 10^{4+4k} \log \log X.
\end{array}
\right .
\end{split}
\end{align}
  The assertion of the proposition now follows from \eqref{equ:bd-S-2}--\eqref{M1bound2}. 
\end{proof}

   Now, the above proposition allows us to establish the following result concerning upper bounds for shifted moments of angular Hecke $L$-functions.
\begin{lemma}
\label{prop: upperboundquad}
With the notation as above and assuming the truth of GRH for $L(s,\lambda_{M_Kl})$ for any $1 \leq l \leq X$. Let $X>0$ be large and let $k\geq 1$ be a fixed integer. Let ${\bf a}=(a_1,\ldots ,a_{k}),\ t=(t_1,\ldots ,t_{k})$ be real $k$-tuples such that $a_i \geq 0$ for all $i$ and that $|t_i|\leq X^A$ for some constant $A>0$.  Then we have
\begin{align*}
   \sum_{\substack{ 1 \leq l \leq X}}\big| L\big(\frac 12+it_1, \lambda_{M_Kl} \big) \big|^{a_1} \cdots \big| L\big(\frac 12+it_{k}, \lambda_{M_Kl} \big) \big|^{a_{k}} \ll_{{\bf a}} &  X(\log X)^{O(1)}.
\end{align*}
\end{lemma}
\begin{proof}
  We first note that by H\"older's inequality, it suffice to show that for fixed $t \in \mr$ with $|t| \leq X^A$ and any real $k > 0$, 
\begin{align*}
   \sum_{\substack{ 1 \leq l \leq X}}\big| L\big(\frac 12+it, \lambda_{M_Kl} \big) \big|^{2k} \ll_{k,A} &  X(\log X)^{O(1)}.
\end{align*}  
   
   Integration by parts implies that
\begin{align}
\label{momentint}
  \sum_{\substack{ 1 \leq l \leq X}}\big| L\big(\frac 12+it, \lambda_{M_Kl} \big) \big|^{k}  =& -\int\limits_{-\infty}^{+\infty}\exp (2kV) \dif \mathcal{N}_t(V,X)
  = 2k\int\limits_{-\infty}^{+\infty}\operatorname{exp}(2kV)\mathcal{N}_t(V,X) \dif V. 
\end{align}
 As $N_t(V, X) \ll X$, we see that
\begin{align*}
 2k\int\limits_{-\infty}^{10\sqrt{\log \log X}}\operatorname{exp}(2kV)\mathcal{N}_t(V,X) \dif V \ll X\int\limits_{-\infty}^{10\sqrt{\log \log X}}\operatorname{exp}(2kV) \dif V \ll X(\log X)^{k^2}.
\end{align*}
 
  On the other hand, upon setting $x = \log X$ in \eqref{logLupperboundnonquad} and bounding the sums  over $\varpi$ there trivially, we see that $\mathcal{N}_t(V,X) =0$  for $V > 10(A+1)\log X/\log \log X$.
It therefore remains to consider the $V$-range with $10\sqrt{\log \log X} <  V \leq  10(A+1)\log X/\log \log X$. We deduce from Proposition \ref{propNbound} that
\begin{align}
 \label{equ:rough-01}
  \mathcal{N}_t(V, X)\ll
\begin{cases}
  X(\operatorname{log}X)^{o(1)}\operatorname{exp}\left(-\frac{V^2}{\log \log X}\right), \quad 10\sqrt{\log \log X} \leq V \leq 10^{4+4k} \log \log X, \\ \\
   X(\operatorname{log}X)^{o(1)}\operatorname{exp}(-(2+4k)V), \quad  10^{4+4k} \log \log X< V \leq  10(A+1)\log X/\log \log X.
\end{cases}
\end{align}
   Applying the bounds given in \eqref{equ:rough-01} to evaluate the integral in \eqref{momentint}, we see that the assertion of the lemma follows. 
\end{proof}

\section{Proof of Theorem \ref{t1}}

   Our proof follows the approach of that of \cite[Theorem 1]{Szab}. Based on the method of A. J. Harper \cite{Harper}, we define for a large number $M$, depending on ${\bf a}$ only, such that
\begin{align*}
\begin{split} 
 \alpha_{0} = 0, \;\;\;\;\; \alpha_{i} = \frac{20^{i-1}}{(\log\log X )^{2}} \;\;\; \mbox{for} \; i \geq 1 \quad \mbox{and} \quad
\mathcal{J} = \mathcal{J}_{{\bf a}, A, X} = 1 + \max\left\{i : \alpha_{i} \leq 10^{-M} \right\} . 
\end{split}
\end{align*}

  We set $x=X^{\alpha_j}$ in Corollary \ref{lemmasum} to see that for $1 \leq l \leq X$,
\begin{align}
\label{Ldecomp}
\begin{split}
 & \sum^{k}_{m=1}a_m\log |L(\tfrac 12+it_m, \lambda_{M_Kl} )| \le 2 \Re \sum^{j}_{i=1} {\mathcal M}_{i,j}(l)+\sum_{0 \leq u \leq \frac{\log\log X}{\log 2}}\Re P_{u}(l )+2(A+1)a\alpha^{-1}_j+O(1),
\end{split}
\end{align}
  where 
\begin{align*}
\begin{split}
  {\mathcal M}_{i,j}(l) = & \sum_{\substack{ (\varpi) \\ X^{\alpha_{i-1}} < N(\varpi) \leq X^{\alpha_{i}}}}  \frac{h(\varpi)\lambda_{M_Kl}(\varpi)}{N(\varpi)^{1/2+1/\log X^{\alpha_{j}}}} \frac{\log (X^{\alpha_{j}}/N(\varpi))}{\log X^{\alpha_{j}}},  \ 1\leq i \leq j \leq \mathcal{J}, \\
   P_{u}(l )= & \sum_{\substack{ (\varpi) \\2^{u} < N(\varpi) \leq 2^{u+1}}}
 \frac{h(\varpi^2)\lambda_{M_Kl}^2(\varpi)}{N(\varpi)}, \quad  0 \leq u \leq \frac{\log\log X}{\log 2}. 
\end{split}
\end{align*}

  We further define the following sets. 
\begin{align*}
  \mathcal{S}(0) =& \left\{1 \leq l \leq X : | \Re {\mathcal M}_{1,r}(l)| > \alpha_{1}^{-3/4} \; \; \mbox{for some}  \; 1 \leq r  \leq \mathcal{J} \right\}, \\
 \mathcal{S}(j) =& \left\{  1 \leq l \leq X  : | \Re {\mathcal M}_{i,r}(l)| \leq \alpha_{i}^{-3/4} \; \; \mbox{for all}  \; 1 \leq i \leq j, \; \mbox{and} \; i \leq r \leq \mathcal{J}, \right. \\
 & \hspace*{4cm} \left. \;\;\;\;\; \text{but }  |\Re{\mathcal M}_{j+1,r}(l)| > \alpha_{j+1}^{-3/4} \; \text{ for some } j+1 \leq r \leq \mathcal{J} \right\}, 
 \quad 1\leq j < \mathcal{J},  \\
 \mathcal{S}(\mathcal{J}) =& \left\{ 1 \leq l \leq X  : |\Re{\mathcal M}_{i, \mathcal{J}}(\chi)| \leq \alpha_{i}^{-3/4} \; \mbox{for all}  \; 1 \leq i \leq
 \mathcal{J} \right\}, \\
 \mathcal{P}(u) =&  \left\{ 1 \leq l \leq X  : |\Re P_{u}(l)| > 2^{-u/10} , \; \text{but} \; |\Re P_{n}(l)| \leq 2^{-n/10} \; \mbox{for all} \; u+1 \leq n \leq \frac{\log\log X}{\log 2} \right\}.
\end{align*}

   Since the treatment for the case $l \not \in \mathcal{P}(u)$ for any $u$ is easier compared to the other cases, we may assume that $l \in \mathcal{P}(u)$ for some 
   $u$.   Observe that
$$ \mathcal{P}(u)= \bigcup_{u=0}^{\log \log X/2}\bigcup_{j=0}^{ \mathcal{J}} \Big (\mathcal{S}(j)\bigcap \mathcal{P}(u) \Big ).$$
 Thus it suffices to prove that
\begin{align}
\label{sumovermj}
\begin{split}
  & \sum_{u=0}^{\log \log X/2}\sum_{j=0}^{\mathcal{J}}\sum_{l \in \mathcal{S}(j)\bigcap \mathcal{P}(u)} \big| L\big( \tfrac12+it_1, \lambda_{M_Kl} \big ) \big|^{a_1} \cdots \big| L\big( \tfrac12+it_{k}, \lambda_{M_Kl}  \big) \big|^{a_{k}} \\
   \ll & X(\log X)^{3(a_1^2+\cdots +a_{k}^2)/8} \\
  & \times \prod_{1\leq j<m \leq k} \Big|\zeta \Big(1+i(t_j-t_m)+\frac 1{\log X} \Big) \Big|^{3a_ja_m/4}\Big|\zeta \Big(1+i(t_j+t_m)+\frac 1{\log X} \Big) \Big|^{a_ja_m/2} \\
  & \times \prod_{1\leq j\leq k} \Big|\zeta \Big(1+2it_j+\frac 1{\log X} \Big) \Big|^{a^2_j/8+a_j/2} \\
  & \times \prod_{1\leq j<m \leq k} \Big|L \Big(1+i(t_j-t_m)+\frac 1{\log X}, \chi^{(D_K)} \Big) \Big|^{3a_ja_m/4}\Big|L\Big(1+i(t_j+t_m)+\frac 1{\log X}, \chi^{(D_K)} \Big) \Big|^{a_ja_m/2} \\
  & \times \prod_{1\leq j\leq k} \Big|L \Big(1+2it_j+\frac 1{\log X}, \chi^{(D_K)} \Big) \Big|^{a^2_j/8-a_j/2}. 
\end{split}
\end{align}

   Observe that
\begin{align*}
\text{meas}(\mathcal{P}(u)) \leq \sum_{\substack{1 \leq l \leq X}}
\Big (2^{u/10} |P_u(l)| \Big )^{2\lceil 2^{u/2}\rceil }.
\end{align*}
We follow the same approach in the proof of Lemma \ref{lem:2.5} and the estimation \eqref{binombound} to bound the right-hand side of the above expression
to see that for $u \geq 10$,
\begin{align}
\label{Pmest}
\begin{split}
  \text{meas}(\mathcal{P}(u)) \ll & X \Big (2^{u/10} \Big )^{2\lceil 2^{u/2}\rceil }\sum^{\lceil 2^{u/2}\rceil}_{i=0}\binom {\lceil 2^{u/2}\rceil}{i}^2 i! \left(\sum_{\substack{(\varpi) \\ 2^u \leq N(\varpi) }} \frac {a^2}{N^2(\varpi)}\right)^{2i} \left(\sum_{2^u \leq p } \frac {a}{p^2}\right)^{2\lceil 2^{u/2}\rceil -2i}\\
  & +\Big (2^{u/10} \Big )^{2\lceil 2^{u/2}\rceil } \cdot (2^{u+1})^{2\lceil 2^{u/2}\rceil} \left( \sum_{\substack{(\varpi) \\ 2^u \leq N(\varpi)}} \frac {a^2}{N^2(\varpi)} \right)^{\lceil 2^{u/2}\rceil } \\
  \ll &  X \Big (2^{u/10} \Big )^{2\lceil 2^{u/2}\rceil }\sum^{\lceil 2^{u/2}\rceil}_{i=0}\binom {\lceil 2^{u/2}\rceil}{i}^2 i! (a^2 2^{-u})^{2\lceil 2^{u/2}\rceil},
\end{split}
 \end{align}
   where the last estimation above follows by noting that $ (2^{u+1})^{2\lceil 2^{u/2}\rceil}(2^{u}/a^2)^{\lceil 2^{u/2}\rceil} \ll X$. 
  
   We deduce from \eqref{Stirling} and \eqref{Pmest}  that
\begin{align}
\label{Pmest1}
\begin{split}
  \text{meas}(\mathcal{P}(u)) 
  \ll &  X \Big (2^{u/10} \Big )^{2\lceil 2^{u/2}\rceil }\sum^{\lceil 2^{u/2}\rceil}_{i=0}\left( \frac{e\lceil 2^{u/2}\rceil}{i} \right)^{2i} i \left( \frac ie \right)^i (a^2 2^{-u})^{2\lceil 2^{u/2}\rceil} \\
  \ll &  X \Big (2^{u/10} \Big )^{2\lceil 2^{u/2}\rceil }\cdot \lceil 2^{u/2}\rceil^2 \cdot (e\lceil 2^{u/2}\rceil)^{2\lceil 2^{u/2}\rceil}  (a^2 2^{-u})^{2\lceil 2^{u/2}\rceil} \\
  \ll & X  2^{-2^{u/2}}.
\end{split}
 \end{align}
 
  Observe that by \cite[Theorem 6.7]{MVa1}, we have for $|t_j|, |t_l|\leq X^A$ and $X$ large enough, 
\begin{align}
\label{zetalowerbound}
 \big|\zeta(1+i(t_j-t_l)+\frac 1{\log X}) \big| \gg \min (\log X, \frac {1}{\log (|t_j-t_l|+4)}) \gg \frac 1{\log X}.
\end{align}  
  
   Moreover, we note that $\chi^{(D_K)}(-1)=-1$ (see Section 9.3 of \cite{MVa1}) so that by a result of M. E. Low \cite{Low68} (see also \cite{Watkins04}) 
   that there is no exceptional zero for $L(s,\chi^{(D_K)})$. This allows us to apply \cite[Theorem 11.4]{MVa1} to see that for $|t_j|, |t_l|\leq X^A$ and $X$ large enough, 
\begin{align}
\label{Llowerbound}
 1=\max (1, \frac {1}{\log (|t_j-t_l|+4)}) \gg \big|\zeta(1+i(t_j-t_l)+\frac 1{\log X}, \chi^{(D_K)}) \big| \gg \min (1, \frac {1}{\log (|t_j-t_l|+4)}) \gg \frac 1{\log X}.
\end{align}  
   
    We then apply the Cauchy-Schwarz inequality, the estimations given in \eqref{Pmest1}--\eqref{Llowerbound} and Lemma \ref{prop: upperboundquad} to see that when $2^{u} \geq (\log\log X)^{3}$,
\begin{align*}
& \sum_{l \in  \mathcal{P}(u) }\big| L\big(\tfrac12+it_1, \lambda_{M_Kl} \big) \big|^{a_1} \cdots \big| L\big(\tfrac12+it_{k}, \lambda_{M_Kl} \big) \big|^{a_{k}} \\
\leq & \left( \text{meas}(\mathcal{P}(u)) \cdot
\sum_{1 \leq l \leq X} \big| L\big(\tfrac12+it_1, \lambda_{M_Kl} \big) \big|^{2a_1} \cdots \big| L\big(\tfrac12+it_{k}, \lambda_{M_Kl} \big) \big|^{2a_{k}} \right)^{1/2}
 \\
 \ll & \left( X \exp\left( -(\log 2)(\log\log X)^{3/2} \right) X (\log X)^{O(1)} \right)^{1/2} \\
  \ll & X(\log X)^{3(a_1^2+\cdots +a_{k}^2)/8} \\
  & \times \prod_{1\leq j<m \leq k} \Big|\zeta \Big(1+i(t_j-t_m)+\frac 1{\log X} \Big) \Big|^{3a_ja_m/4}\Big|\zeta \Big(1+i(t_j+t_m)+\frac 1{\log X} \Big) \Big|^{a_ja_m/2} \\
  & \times \prod_{1\leq j\leq k} \Big|\zeta \Big(1+2it_j+\frac 1{\log X} \Big) \Big|^{a^2_j/8+a_j/2} \\
  & \times \prod_{1\leq j<m \leq k} \Big|L \Big(1+i(t_j-t_m)+\frac 1{\log X}, \chi^{(D_K)} \Big) \Big|^{3a_ja_m/4}\Big|L\Big(1+i(t_j+t_m)+\frac 1{\log X}, \chi^{(D_K)} \Big) \Big|^{a_ja_m/2} \\
  & \times \prod_{1\leq j\leq k} \Big|L \Big(1+2it_j+\frac 1{\log X}, \chi^{(D_K)} \Big) \Big|^{a^2_j/8-a_j/2}. 
\end{align*}

   The above implies that we may further assume that $0 \leq u \leq (3/\log 2)\log\log\log X$.
We next note that, 
\begin{align}
\label{S0est}
\begin{split}
\text{meas}(\mathcal{S}(0)) \ll & \sum_{1 \leq l \leq X} 
\sum^{\mathcal{J}}_{j=1}
\Big ( \alpha^{3/4}_{1}{|\mathcal
M}_{1, j}(l)| \Big)^{2\lceil 1/(10\alpha_{1})\rceil } =\sum^{\mathcal{J}}_{j=1}\sum_{1 \leq l \leq X} \Big ( \alpha^{3/4}_{1}{|\mathcal
M}_{1, j}(l)| \Big)^{2\lceil 1/(10\alpha_{1})\rceil }.
\end{split}
\end{align}

  Denote $w_K(n)$ the multiplicative function such that $w_K(\varpi^{\alpha}) = \alpha!$ for prime powers $\varpi^{\alpha}$ and $w_K(u)=1$ for $u \in U_K$. 
  Denote also $\Omega_K(n)$ the number of distinct prime ideal powers dividing $(n)$. We have for fixed $1 \leq j \leq X$, 
\begin{align}
\label{suml}
\begin{split}
& \sum_{1 \leq l \leq X} \Big ( \alpha^{3/4}_{1}{|\mathcal
M}_{1, j}(l)| \Big)^{2\lceil 1/(10\alpha_{1})\rceil } \\
=& ( \alpha^{3/4}_{1} )^{2\lceil 1/(10\alpha_{1})\rceil }(\lceil 1/(10\alpha_{1})\rceil !)^2 \sum_{1 \leq l \leq X}\sum_{ \substack{ (n_i), (n'_i)\\ \Omega_K(n_i)=\Omega_K(n'_i) = \lceil 1/(10\alpha_{1})\rceil \\ \varpi|n_in'_i \implies  X^{\alpha_{0}} < N(\varpi) \leq X^{\alpha_{1}} }} \frac{H_j(n_i)\overline{H_j(n'_i)}\lambda_{M_Kl}(n_i\overline{n'_i}) }{\sqrt{N(n_in'_i)} w(n_i) w(n'_i)},
\end{split}
\end{align}  
  where $H_j(n)$ is a totally multiplicative function of $n$ such that $H_j(u)=1$ for $u \in U_K$ and 
\begin{align}
\label{Hjdef}
\begin{split}
 H_j(\varpi)=\frac {h(\varpi)\log (X^{\alpha_j}/N(\varpi))}{N(\varpi)^{1/\log X^{\alpha_{j}}}\log X^{\alpha_j}}. 
\end{split}
\end{align}    

  Note that we have $|H_j(\varpi)| \leq h(\varpi) \leq a$. We now apply Lemma \ref{prsum} to evaluate the right-hand side expression of \eqref{suml} to see that
\begin{align}
\label{suml1}
\begin{split}
& \sum_{1 \leq l \leq X} \Big ( \alpha^{3/4}_{1}{|\mathcal
M}_{1, j}(l)| \Big)^{2\lceil 1/(10\alpha_{1})\rceil } \\
\ll & ( a\alpha^{3/4}_{1} )^{2\lceil 1/(10\alpha_{1})\rceil }(\lceil 1/(10\alpha_{1})\rceil !)^2 X\sum_{ \substack{ (n_i), (n'_i)\\ \Omega_K(n_i)=\Omega_K(n'_i) = \lceil 1/(10\alpha_{1})\rceil \\ \varpi|n_in'_i \implies  X^{\alpha_{0}} < N(\varpi) \leq X^{\alpha_{1}} \\ n_i\overline{n'_i} \in \mn }} \frac{1}{\sqrt{N(n_in'_i)} w_K(n_i) w_K(n'_i)} \\
& +  ( a\alpha^{3/4}_{1} )^{2\lceil 1/(10\alpha_{1})\rceil }(\lceil 1/(10\alpha_{1})\rceil !)^2 \sum_{ \substack{ (n_i), (n'_i)\\ \Omega_K(n_i)=\Omega_K(n'_i) = \lceil 1/(10\alpha_{1})\rceil \\ \varpi|n_in'_i \implies  X^{\alpha_{0}} < N(\varpi) \leq X^{\alpha_{1}} }} \frac{1}{w_K(n_i) w_K(n'_i)}.
\end{split}
\end{align}    
  
  It follows from the above and \eqref{Stirling} that 
\begin{align}
\label{suml1-1}
\begin{split}
 & ( a\alpha^{3/4}_{1} )^{2\lceil 1/(10\alpha_{1})\rceil }(\lceil 1/(10\alpha_{1})\rceil !)^2 \sum_{ \substack{ (n_i), (n'_i)\\ \Omega_K(n_i)=\Omega_K(n'_i) = \lceil 1/(10\alpha_{1})\rceil \\ \varpi|n_in'_i \implies  X^{\alpha_{0}} < N(\varpi) \leq X^{\alpha_{1}} }} \frac{1}{w_K(n_i) w_K(n'_i)} \\
 =& ( a\alpha^{3/4}_{1} )^{2\lceil 1/(10\alpha_{1})\rceil }\lceil 1/(10\alpha_{1})\rceil (\frac {\lceil 1/(10\alpha_{1})\rceil}{e})^{\lceil 1/(10\alpha_{1})\rceil}(\sum_{ \substack{ (\varpi)\\ X^{\alpha_{0}} < N(\varpi) \leq X^{\alpha_{1}} }}1)^{2\lceil 1/(10\alpha_{1})\rceil } \\
 \ll & X^{1/4}. 
\end{split}
\end{align}    

   Note also that, when $n_i\overline{n'_i} \in \mn$,  then for any prime $\varpi$ satisfying $N(\varpi)=p$ a rational prime and $\varpi^k \| n_i\overline{n'_i}$ 
   for some integer $k \geq 1$, we must also have  $\overline{\varpi}^k \| n_i\overline{n'_i}$, so that $p^{2k} \mid N(n_i\overline{n'_i})$. We write for this case that
   $n_i=\varpi^h\overline{\varpi}^gm_i$ with $(n_i, m_i)=1$, then we must have $n'_i=\varpi^{k-h}\overline{\varpi}^{k-g}m'_i$ with $(n'_i, m'_i)=1$ and 
   therefore the expression $(\sqrt{N(n_in'_i)} w(n_i) w(n'_i))^{-1}$ has a factor 
\begin{align*}
\begin{split}
 p^{-k} \sum^k_{h=0}\sum^k_{g=0}\frac {1}{h!(k-h)!g!(k-g)!}= p^{-k}\frac {1}{(k!)^2}\big (\sum^{k}_{i=1}\binom {k}{i}\big)^2=p^{-k}\frac {2^k}{(k!)^2}
 \leq 2p^{-k}\sum^k_{i=0}\frac {1}{i!(k-i)!}
\end{split}
\end{align*}

  Note also that we have $N(p)=p^2$ for a rational prime $p$. It follows that we have
\begin{align}
\label{firsttermest}
\begin{split}
& (\lceil 1/(10\alpha_{1})\rceil !)^2\sum_{ \substack{ (n_i), (n'_i)\\ \Omega_K(n_i)=\Omega_K(n'_i) = \lceil 1/(10\alpha_{1})\rceil\\ \varpi|n_in'_i \implies  X^{\alpha_{0}} < N(\varpi) \leq X^{\alpha_{1}} \\ n_i\overline{n'_i} \in \mn }} \frac{1}{\sqrt{N(n_in'_i)} w_K(n_i) w_K(n'_i)} \\
\leq & \left(\sum_{\substack{ p \\ X^{\alpha_{0}/2} < p \leq X^{\alpha_{1}/2} }}  \frac {2}{p} \right)^{ \lceil 1/(10\alpha_{1})\rceil}\cdot 
\left(\sum_{\substack{ p \\ X^{\alpha_{0}/2} < p \leq X^{\alpha_{1}/2} }}  \frac {2}{p} \right)^{ \lceil 1/(10\alpha_{1})\rceil} \\
\leq & \left(2\log \log X+O(1)\right)^{ 2\lceil 1/(10\alpha_{1})\rceil},
\end{split}
\end{align}
  where the last estimation above follows from \eqref{merten}. 
  
  Observe that we have 
\begin{align}
\label{sump1}
 \mathcal{J} \leq \log\log\log X , \; \alpha_{1} = \frac{1}{(\log\log X)^{2}} , \; \mbox{and} \; \sum_{\substack{(\varpi) \\ N(\varpi) \leq X^{1/(\log\log X)^{2}}}} \frac{1}{N(\varpi)} \leq \log\log X ,
\end{align}
  where the last estimation also follows from \eqref{merten}.   
  
  We derive from \eqref{firsttermest} and \eqref{sump1} that 
\begin{align}
\label{firsttermest1}
\begin{split}
 ( a\alpha^{3/4}_{1} )^{2\lceil 1/(10\alpha_{1})\rceil }(\lceil 1/(10\alpha_{1})\rceil !)^2 X\sum_{ \substack{ (n_i), (n'_i)\\ \Omega_K(n_i)=\Omega_K(n'_i) = \lceil 1/(10\alpha_{1})\rceil \\ \varpi|n_in'_i \implies  X^{\alpha_{0}} < N(\varpi) \leq X^{\alpha_{1}} \\ n_i\overline{n'_i} \in \mn }} \frac{1}{\sqrt{N(n_in'_i)} w_K(n_i) w_K(n'_i)} \ll X e^{-(\log\log X)^{2}/10}  .
\end{split}
\end{align}    

  Applying the estimations obtained in \eqref{suml1-1} and \eqref{firsttermest1} in \eqref{suml1} and then substituting the result in \eqref{S0est}, we see that 
\begin{align*}
\text{meas}(\mathcal{S}(0)) \ll &
\mathcal{J}X e^{-1/\alpha_{1}}\ll X e^{-(\log\log X)^{2}/10}  .
\end{align*}
  We then apply the Cauchy-Schwarz inequality and Lemma \ref{prop: upperboundquad} to deduce that
\begin{align*}
& \sum_{l \in  \mathcal{S}(0)} \big| L\big(\tfrac12+it_1, \lambda_{M_Kl} \big) \big|^{a_1} \cdots \big| L\big(\tfrac12+it_{k}, \lambda_{M_Kl} \big) \big|^{a_{k}} \\
 \leq &   \left( \text{meas}(\mathcal{S}(0)) \cdot
\sum_{\substack{1 \leq l \leq X}}  \big| L\big(\tfrac12+it_1, \lambda_{M_Kl} \big) \big|^{2a_1} \cdots \big| L\big(\tfrac12+it_{k}, \lambda_{M_Kl} \big) \big|^{2a_{k}} \right)^{1/2}
 \\
 \ll & \left( X \exp\left( -(\log\log X)^{2}/10 \right) X (\log X)^{O(1)} \right)^{1/2} \ll X(\log X)^{k^2}.
\end{align*}

 The above implies that we may further assume that $j \geq 1$. We deduce from \eqref{Ldecomp} that for $l \in \mathcal{S}(j)$, 
\begin{align*}
\begin{split}
 & \big|L\big(\tfrac12+it_1, \lambda_{M_Kl} \big) \big|^{a_1} \cdots \big| L\big(\tfrac12+it_{k}, \lambda_{M_Kl} \big) \big|^{a_{k}} \ll \exp \left(\frac {2(A+1)a}{\alpha_j} \right) \exp \Big (
 2\Re\sum^j_{i=0}{\mathcal M}_{i,j}(l)+\Re\sum^{\log \log X/\log 2}_{u=0}P_u(l) \Big ).
\end{split}
 \end{align*}

   When estimating the sum of $\big|L\big(\tfrac12+it_1, \lambda_{M_Kl} \big) \big|^{a_1} \cdots \big| L\big(\tfrac12+it_{k}, \lambda_{M_Kl} \big) \big|^{a_{k}}$
    over $\mathcal{S}(j)\bigcap \mathcal{P}(u)$, we need to separate the sums over $\varpi$ with $N(\varpi) \leq 2^{u+1}$ on the right-hand side of the above 
    expression from those with $N(\varpi)>2^{u+1}$. For this, we observe that when $l \in \mathcal{P}(u)$, we have
\begin{align*}
\begin{split}
 & 2\Re \sum_{\substack{ (\varpi) \\ N(\varpi) \leq 2^{u+1}}}  \frac{h(\varpi)\lambda_{M_Kl}(\varpi)}{N(\varpi)^{1/2+1/\log X^{\alpha_{j}}}} \frac{\log (X^{\alpha_{j}}/N(\varpi))}{\log X^{\alpha_{j}}} +
  \Re \sum_{\substack{ (\varpi) \\ N(\varpi) \leq \log X}}
 \frac{h(\varpi^2)\lambda_{M_Kl}^2(\varpi)}{N(\varpi)}
   \\
 \leq & 2\Re\sum_{\substack{ (\varpi) \\ N(\varpi) \leq 2^{u+1}}}  \frac{h(\varpi)\lambda_{M_Kl}(\varpi)}{N(\varpi)^{1/2+1/\log X^{\alpha_{j}}}} \frac{\log (X^{\alpha_{j}}/N(\varpi))}{\log X^{\alpha_{j}}} +
 \Re \sum_{\substack{ (\varpi) \\ N(\varpi) \leq 2^{u+1}}}
 \frac{h(\varpi^2)\lambda_{M_Kl}^2(\varpi)}{N(\varpi)}+O(1) \leq a2^{u/2+3}+O(1).
\end{split}
 \end{align*}

 We deduce from the above that
\begin{align}
\label{LboundinSP}
\begin{split}
   &\sum_{l \in \mathcal{S}(j)\bigcap \mathcal{P}(u)} \big|L\big(\tfrac12+it_1, \lambda_{M_Kl} \big) \big|^{a_1} \cdots \big| L\big(\tfrac12+it_{k}, \lambda_{M_Kl} \big) \big|^{a_{k}} \\
   \ll & e^{a2^{u/2+4}} \sum_{l \in \mathcal{S}(j)\bigcap \mathcal{P}(u)}
 \exp \left(\frac {2(A+1)a}{\alpha_j} \right)\exp \Big ( 2\Re \sum_{\substack{ (\varpi) \\ 2^{u+1}<N(\varpi) \leq X^{\alpha_j}}}
 \frac{h(\varpi)\lambda_{M_Kl}(\varpi)}{N(\varpi)^{1/2+1/\log X^{\alpha_{j}}}} \frac{\log (X^{\alpha_{j}}/N(\varpi))}{\log X^{\alpha_{j}}}\Big )
  \\
\ll &  e^{a2^{u/2+4}} \exp \left(\frac {2(A+1)a}{\alpha_j} \right)\sum_{l \in \mathcal{S}(j)} \Big (2^{u/10}|P_u(\chi)| \Big )^{2\lceil 2^{u/2}\rceil }
\exp \Big ( 2\Re{\mathcal M}'_{1,j}(l)+2\Re \sum^j_{i=2}{\mathcal M}_{i,j}(l)\Big ),
\end{split}
 \end{align}
   where 
\begin{align*}
\begin{split}
 {\mathcal M}'_{1,j}(l) :=  \sum_{\substack{ (\varpi) \\ 2^{u+1}<N(\varpi) \leq X^{\alpha_1}}}
  \frac{h(\varpi)\lambda_{M_Kl}(\varpi)}{N(\varpi)^{1/2+1/\log X^{\alpha_{j}}}} \frac{\log (X^{\alpha_{j}}/N(\varpi))}{\log X^{\alpha_{j}}}.
\end{split}
 \end{align*}

    Note that for $0 \leq u \leq (3/\log 2)\log\log\log X$ and $X$ large enough, 
\begin{align*}
\begin{split}
  \sum_{\substack{ (\varpi) \\ N(\varpi) \leq 2^{u+1}}}
  \frac{h(\varpi)\lambda_{M_Kl}(\varpi)}{N(\varpi)^{1/2+1/\log X^{\alpha_{j}}}} \frac{\log (X^{\alpha_{j}}/N(\varpi))}{\log X^{\alpha_{j}}}\leq \sum_{\substack{ (\varpi) \\ N(\varpi) \leq 2^{u+1}}}
 \frac{a}{\sqrt{N(\varpi)}} \leq \frac {100a \cdot 2^{u/2}}{u+1} \leq 100a(\log \log X)^{3/2}(\log \log \log X)^{-1},
\end{split}
 \end{align*}
   where the last estimation above follows from \eqref{PIT} and partial summation. 

  The above implies that when $l \in \mathcal{S}(j)$ and $X$ large enough,
\begin{align*}
\begin{split}
 |{\mathcal M}'_{1,j}(l)|\leq 100a(\log \log X)^{3/2}(\log \log \log X)^{-1}+{\mathcal M}_{1,j}(\chi) \leq 1.01\alpha^{-3/4}_1=1.01(\log \log X)^{3/2}.
\end{split}
 \end{align*}

  As we also have $|{\mathcal M}_{i, j}(l)| \leq  \alpha^{-3/4}_i$ when $l \in \mathcal{S}(j)$, we may therefore apply \cite[Lemma 5.2]{Kirila} to see that
\begin{align*}
\begin{split}
\exp \Big (2\Re {\mathcal M}'_{1,j}(l)+2\Re\sum^j_{i=2}{\mathcal M}_{i,j}(l)\Big ) \ll \Big | E_{e^2\alpha^{-3/4}_1}\big({\mathcal M}'_{1,j}(l)\big) \Big |^2
\prod^j_{i=2}\Big | E_{e^2\alpha^{-3/4}_i}\big({\mathcal M}_{i,j}(l)\big)\Big |^2,
\end{split}
 \end{align*}
  where 
\begin{align*}
  E_{\ell}(x) = \sum_{j=0}^{\lceil\ell \rceil} \frac {x^{j}}{j!}
\end{align*}

   We then deduce from \eqref{LboundinSP} and the description on $\mathcal{S}(j)$ that when $j \geq 1$,
\begin{align}
\label{L2k}
\begin{split}
  & \sum_{l \in \mathcal{S}(j)\bigcap \mathcal{P}(u)} \big|L\big(\tfrac12+it_1, \lambda_{M_Kl} \big) \big|^{a_1} \cdots \big| L\big(\tfrac12+it_{k}, \lambda_{M_Kl} \big) \big|^{a_{k}}  \\
 \ll & e^{a2^{u/2+4}}\exp \left( \frac {2(A+1)a}{\alpha_j} \right)
 \sum^{\mathcal{J} }_{r=j+1} \sum_{\substack{1 \leq l \leq X}} \Big (2^{u/10}|P_u(l)| \Big )
 ^{2\lceil 2^{u/2}\rceil } \\
& \hspace*{2cm} \times \exp \Big (  2\Re {\mathcal M}'_{1,j}(l)+2\Re \sum^j_{i=2}{\mathcal M}_{i,j}(l)\Big )\Big ( \alpha^{3/4}_{j+1}|{\mathcal
M}_{j+1, r}(l)|\Big)^{2\lceil 1/(10\alpha_{j+1})\rceil } \\
\ll &  e^{a2^{u/2+4}}\exp \left(\frac {2(A+1)a}{\alpha_j} \right) \sum^{\mathcal{J} }_{r=j+1}\mathcal S_r. 
\end{split}
 \end{align}
  where
\begin{align}
\label{Srdef}
\begin{split}
 \mathcal S_r :=& \sum_{\substack{1 \leq l \leq X}} 
\Big (2^{u/10}|P_u(l)| \Big )^{2\lceil 2^{u/2}\rceil }\Big | E_{e^2\alpha^{-3/4}_1}\big({\mathcal M}'_{1,j}(l)\big) \Big |^2
\prod^j_{i=2}\Big | E_{e^2\alpha^{-3/4}_i}\big({\mathcal M}_{i,j}(l)\big) \Big |^2\Big ( \alpha^{3/4}_{j+1}|{\mathcal
M}_{j+1, r}(l)|\Big)^{2\lceil 1/(10\alpha_{j+1})\rceil } .
\end{split}
\end{align}

   We now evaluate of $S_r$ for a fixed $r$ by expanding the factors in $S_r$ into Dirichlet series.  For this, we define functions $b_j(n)$, $1 \leq j \leq \mathcal{J}$ such that $b_j(n)=0$ or $1$ and $b_j(n)=1$ if and only if $\Omega_K(n) \leq e^2\alpha^{-3/4}_j$ and the primes $\varpi$ dividing $n$ satisfying $X^{\alpha_{j-1}}<N(\varpi) \leq X^{\alpha_j}$ for $j \geq 2$ and $2^{u+1}< N(\varpi) \leq X^{\alpha_1}$ for $j=1$. 
    We use these notations to see that
\begin{equation*}
E_{e^2\alpha^{-3/4}_i}( {\mathcal M}_{i,j}(l)) = \sum_{ \substack {(n_i)}} \frac{H_j(n_i)}{\sqrt{N(n_i)}w(n_i)}  b_i(n_i)\lambda_{M_Kl}(n_i) , \quad 2\le i\le j.
\end{equation*}
  Note that a similar expression holds for $E_{e^2\alpha^{-3/4}_1}( {\mathcal M}'_{1,j}(l))$ as well.
  
   Note that each Dirichlet polynomial $E_{e^2\alpha^{-3/4}_i}( {\mathcal M}_{i,j}(l)), 2 \leq i \leq j$ is a short Dirichlet polynomial of length at most $X^{\alpha_{i}\cdot e^2\alpha^{-3/4}_i}=X^{e^2\alpha^{1/4}_i}$. Similarly, the Dirichlet polynomial $E_{e^2\alpha^{-3/4}_1}( {\mathcal M}'_{1,j}(l))$ is a short Dirichlet polynomial of length at most $X^{\alpha_{1}\cdot e^2\alpha^{-3/4}_1}=X^{e^2\alpha^{1/4}_1}$. Note that we have
\begin{align*}
 \sum^{\mathcal{J}}_{i=1} e^2\alpha^{1/4}_i \leq 2e^210^{-M/4}.
\end{align*}   
   The above implies that $E_{e^2\alpha^{-3/4}_1}({\mathcal M}'_{1,j}(l))\prod^{j}_{i=2}E_{e^2\alpha^{-3/4}_i}( {\mathcal M}_{i,j}(l))$ is also a short Dirichlet
    polynomial whose length is at most $X^{2e^210^{-M/4}}$.  We write for simplicity that
\begin{align}
\label{Eprodexpression}
 E_{e^2\alpha^{-3/4}_1}({\mathcal M}'_{1,j}(l))\prod^{j}_{i=2}E_{e^2\alpha^{-3/4}_i}( {\mathcal M}_{i,j}(l))= \sum_{\substack{ (n) \\ N(n)  \leq X^{2e^210^{-M/4}}}} \frac{x_n}{\sqrt{N(n)}} \lambda_{M_Kl}(n),
\end{align}
where $x_n$ denotes the coefficient of the Dirichlet polynomial resulting from the expansion of the left-hand side of \eqref{Eprodexpression}.  Now the observation that $|H_j(n_i)|\leq a^{\Omega_K(n_i)}$ lead to
\begin{align}
\label{xnbounds}
 x_n \ll a^{\sum^{j}_{i=1}e^2\alpha^{-3/4}_i} \ll a^{100e^2(\log \log X)^{3/2}} \ll X^{\varepsilon}.
\end{align}

   On the other hand, 
\begin{align}
\label{Mjuexp}
\begin{split}
 & \Big (   \alpha^{3/4}_{j+1}{\mathcal M}_{j+1,r}(l)\Big )^{\lceil 1/(10\alpha_{j+1})\rceil } \\
= &  \Big (  \alpha^{3/4}_{j+1}\Big)^{\lceil 1/(10\alpha_{j+1})\rceil } \sum_{ \substack{ (n_r) \\ \Omega_K(n_r) = \lceil 1/(10\alpha_{j+1})\rceil \\ \varpi|n_r \implies
X^{\alpha_j}<N(\varpi) \leq X^{\alpha_{j+1}}}}
\frac{(\lceil 1/(10\alpha_{j+1})\rceil  )!H_r(n_r)}{\sqrt{N(n_r)}w(n_r)}\lambda_{M_Kl}(n_r) =:  \sum_{n_r} \frac {y_{n_r}}{\sqrt{N(n_r)}}\lambda_{M_Kl}(n_r).
\end{split}
\end{align}
Note that $N(n_r) \leq (X^{\alpha_{j+1}})^{\lceil 1/(10\alpha_{j+1})\rceil } \leq X^{1/10}$.  The inequality $|H_r(n_r)|\leq a^{\Omega_K(n_r)}$
together with \eqref{Stirling} implies that
\begin{align}
\label{ynbound}
\begin{split}
   y_{n_r} \ll \lceil 1/(10\alpha_{j+1})\rceil \Big ( \frac {a\cdot \alpha^{3/4}_{j+1} \lceil 1/(10\alpha_{j+1})\rceil }{e} \Big)^{\lceil 1 /(10\alpha_{j+1})\rceil } \ll X^{\varepsilon}.
\end{split}
\end{align}

  Similarly, we have
\begin{align}
\label{Puexp}
\begin{split}
 & \Big (2^{u/10}P_u(l) \Big )^{\lceil 2^{u/2}\rceil } \\
= &  \Big ( 2^{u/10} \Big)^{\lceil 2^{u/2}\rceil  } \sum_{ \substack{ (n_u) \\ \Omega_K(n_u) = \lceil 2^{u/2}\rceil \\ \varpi|n_u \implies
2^{u}<N(\varpi) \leq 2^{u+1}}}
\frac{(\lceil 2^{u/2}\rceil  )!G(n_u)}{N(n_u)w(n_u)}\lambda_{M_Kl}(n^2_u) =:  \sum_{n_u} \frac {z_{n_u}}{N(n_u)}\lambda_{M_Kl}(n^2_u),
\end{split}
\end{align}
   where $G(n)$ is a totally multiplicative function of $n$ such that $G(u)=1$ for $u \in U_K$ and 
\begin{align*}
\begin{split}
 G(\varpi)=h(\varpi^2). 
\end{split}
\end{align*}    
  Note that $N(n_u) \leq (2^{u+1})^{\lceil 2^{u/2}\rceil } \leq X^{1/10}$ for $u \leq \log \log X/\log 2$.  The inequality $|G(n_u)|\leq a^{\Omega_K(n_u)}$
together with \eqref{Stirling} implies that
\begin{align}
\label{znbound}
\begin{split}
   z_{n_u} \ll \lceil 2^{u/2}\rceil \Big ( \frac {a\cdot 2^{u/10}\lceil 2^{u/2}\rceil }{e} \Big)^{\lceil 2^{u/2}\rceil } \ll X^{\varepsilon}.
\end{split}
\end{align}

  We deduce from \eqref{Srdef}, \eqref{Eprodexpression}, \eqref{Mjuexp} and \eqref{Puexp} that
\begin{align*}
\begin{split}
 \mathcal S_r =& \sum \sum_{\substack{1 \leq l \leq X}}  \frac {x_n\overline{x_
 {n'}}y_{n_r}\overline{y_{n'_r}}z_{n_u}\overline{z_{n'_u}}}{\sqrt{N(nn'n_rn'_rn^2_u(n'_u)^2)}}
 \lambda_{M_Kl}(n\overline{n'}n_r\overline{n'_r}n^2_u\overline{(n'_u)^2}). 
\end{split}
\end{align*}
  
  We now apply Lemma \ref{prsum} to evaluate the inner sum above. In view of \eqref{xnbounds}, \eqref{ynbound} and \eqref{znbound} 
  and the length of each corresponding Dirichlet polynomials, we see that the main term from \eqref{sumlambda} dominates so that
\begin{align*}
\begin{split}
 \mathcal S_r \ll & X\sum_{\substack{n\overline{n'}n_r\overline{n'_r}n^2_u\overline{(n'_u)^2} \in \mn}}  \frac {x_n\overline{x_{n'}}y_{n_r}\overline{y_{n'_r}}z_{n_u}\overline{z_{n'_u}}}{\sqrt{N(nn'n_rn'_rn^2_u(n'_u)^2)}}.
\end{split}
\end{align*}  
    Note that the elements $nn', n_rn'_r$ and $n^2_u(n'_u)^2$ are mutually co-prime so that $nn'n_rn'_rn^2_u(n'_u)^2 \in \mn$ if and only if $nn' \in \mn, n_rn'_r\in \mn, n^2_u(n'_u)^2\in \mn$. It follows that we have
\begin{align}
\label{Srbounddecomp}
\begin{split}
 \mathcal S_r \ll & X \prod^j_{i=1} \Big (\sum_{\substack{n_i\overline{n'_i} \in \mn}} \frac{H_j(n_i)\overline{H_j(n'_i)}}{\sqrt{N(n_in'_i)}w(n_i)w(n'_i)}  b_i(n_i)b_i(\overline{n'_i}) \Big )   \\
 & \times \Big ( 2^{u/10} \Big)^{2\lceil 2^{u/2}\rceil  } \sum_{ \substack{ (n_u), (n'_u) \\ \Omega_K(n_u) =\Omega_K(n'_u)=  \lceil 2^{u/2}\rceil \\ \varpi|n_un'_u \implies
2^{u}<N(\varpi) \leq 2^{u+1}}}
\frac{((\lceil 2^{u/2}\rceil  )!)^2G(n_u)\overline{G(n'_u)}}{N(n_un'_u)w(n_u)w(n'_u)} \\
& \times  \Big (  \alpha^{3/4}_{j+1}\Big)^{2\lceil 1/(10\alpha_{j+1})\rceil } \sum_{ \substack{ (n_r), (n'_r) \\ \Omega_K(n_r) = \Omega_K(n'_r) = \lceil 1/(10\alpha_{j+1})\rceil \\ \varpi|n_rn'_r \implies
X^{\alpha_j}<N(\varpi) \leq X^{\alpha_{j+1}}}}
\frac{((\lceil 1/(10\alpha_{j+1})\rceil  )!)^2H_r(n_r)\overline{H_r(n'_r)}}{\sqrt{N(n_rn'_r)}w(n_r)w(n'_r)}.
\end{split}
\end{align}      
    
  Similar to \eqref{firsttermest}, we have
\begin{align}
\label{SrbounddecompP}
\begin{split}
 & \Big ( 2^{u/10} \Big)^{2\lceil 2^{u/2}\rceil  } \sum_{ \substack{ (n_u), (n'_u) \\ \Omega_K(n_u) =\Omega_K(n'_u)=  \lceil 2^{u/2}\rceil \\ \varpi|n_un'_u \implies
2^{u}<N(\varpi) \leq 2^{u+1}}}
\frac{((\lceil 2^{u/2}\rceil  )!)^2G(n_u)\overline{G(n'_u)}}{N(n_un'_u)w(n_u)w(n'_u)} \\
\ll & \Big ( a2^{u/10} \Big)^{2\lceil 2^{u/2}\rceil  }\left(\sum_{\substack{ p \\ 2^{u/2} < p \leq 2^{(u+1)/2} }}  \frac {2}{p^2} \right)^{ 2\lceil 2^{u/2}\rceil  }\\
\ll & 2^{- u\lceil 2^{u/2}\rceil/2 }. 
\end{split}
\end{align}     

  Also, 
\begin{align}
\label{SrbounddecompM}
\begin{split}
 &  \Big (  \alpha^{3/4}_{j+1}\Big)^{2\lceil 1/(10\alpha_{j+1})\rceil } \sum_{ \substack{ (n_r), (n'_r) \\ \Omega_K(n_r) = \Omega_K(n'_r) = \lceil 1/(10\alpha_{j+1})\rceil \\ \varpi|n_rn'_r \implies
X^{\alpha_j}<N(\varpi) \leq X^{\alpha_{j+1}}}}
\frac{((\lceil 1/(10\alpha_{j+1})\rceil  )!)^2H_r(n_r)\overline{H_r(n'_r)}}{\sqrt{N(n_rn'_r)}w(n_r)w(n'_r)}\\
\ll & \Big (  a\alpha^{3/4}_{j+1}\Big)^{2\lceil 1/(10\alpha_{j+1})\rceil } \left(\sum_{\substack{ p \\ X^{\alpha_j/2}< p \leq X^{\alpha_{j+1}/2}}}  \frac {2}{p} \right)^{ 2\lceil 1/(10\alpha_{j+1})\rceil  }\\
\ll &  e^{-10^3(A+1)a/(2\alpha_{j+1})}. 
\end{split}
\end{align}       

   It remains to estimate for each fixed $1 \leq i \leq j$,  
\begin{align}
\label{SrbounddecompH}
\begin{split}
 \sum_{\substack{n_i\overline{n'_i} \in \mn}} \frac{H_j(n_i)\overline{H_j(n'_i)}}{\sqrt{N(n_in'_i)}w(n_i)w(n'_i)}  b_i(n_i)b_i(\overline{n'_i})   .
\end{split}
\end{align}  
 
  Note that the factor $ b_i(n_i)$ limits $n_i$ to $\Omega_K(n_i) \leq  e^2\alpha^{-3/4}_i$ and the primes $\varpi$ dividing $n_i$ satisfying $X^{\alpha_{i-1}}<N(\varpi) \leq X^{\alpha_i}$. If we remove the restriction on $\Omega_K(n_i)$ and $\Omega_K(n'_i)$, then \eqref{Srbounddecomp} becomes
\begin{align*}
 \sum_{\substack{n_i\overline{n'_i} \in \mn\\ \varpi | n_in'_i \Rightarrow X^{\alpha_{i-1}}<N(\varpi) \leq X^{\alpha_i}}} \frac{H_j(n_i)\overline{H_j(n'_i)}}{\sqrt{N(n_in'_i)}w(n_i)w(n'_i)}.
\end{align*}
     On the other hand, we apply Rankin's trick to estimate the error introduced in the above process 
     by observing that $2^{\Omega_K(n_i) -e^2\alpha^{-3/4}_i}\ge 1$ when $\Omega_K(n_l) > e^2\alpha^{-3/4}_i$. 
     By the symmetry between $n_i$ and $n'_i$, we see that the error does not exceed
\begin{align*}
\sum_{\substack{n_i\overline{n'_i} \in \mn\\ \varpi | n_in'_i \Rightarrow X^{\alpha_{i-1}}<N(\varpi) \leq X^{\alpha_i}}} \frac{2^{\Omega_K(n_i) -e^2\alpha^{-3/4}_i}|H_j(n_i)\overline{H_j(n'_i)}|}{\sqrt{N(n_in'_i)}w(n_i)w(n'_i)} .
\end{align*}   

    Note that both the main term and the error term  above are now multiplicative functions of $n$ and hence can be evaluated in terms of products over primes.
   Observe also that when $\varpi\overline{\varpi}=p$ is a rational prime, then we can have $\varpi \mid n_i, \overline{\varpi} \mid n'_i$ or $\overline{\varpi} \mid n_i, \varpi \mid n'_i$ or $\varpi\overline{\varpi}|n_i$ or  $\varpi\overline{\varpi}|n'_i$. It follows that
\begin{align*}
\begin{split}
 & \sum_{\substack{n_i\overline{n'_i} \in \mn\\ \varpi | n_in'_i \Rightarrow X^{\alpha_{i-1}}<N(\varpi) \leq X^{\alpha_i}}} \frac{H_j(n_i)\overline{H_j(n'_i)}}{\sqrt{N(n_in'_i)}w(n_i)w(n'_i)} \\
 =& \prod_{\substack{X^{\alpha_{i-1}/2}<p \leq X^{\alpha_i/2} \\ p \ \text{inert} }}\Big  (1+\frac {H_j(p)+\overline{H_j(p)}}{p}+O\Big( \frac 1{p^{2}} \Big) \Big ) \\
 & \times \prod_{\substack{X^{\alpha_{i-1}/2}<p \leq X^{\alpha_i/2} \\ p=\varpi \overline{\varpi} }}\Big  (1+\frac {H_j(\varpi)\overline{H_j(\overline{\varpi}})+H_j(\overline{\varpi})\overline{H_j(\varpi)}+H_j(\varpi)H_j(\overline{\varpi})+\overline{H_j(\varpi)H_j(\overline{\varpi})})}{p}+O\Big( \frac 1{p^{2}} \Big) \Big ).
\end{split}
\end{align*}
 
    Note that $h(\varpi)=h(\overline{\varpi})=h(\sqrt{p})$ when $p=\varpi\overline{\varpi}$ by \eqref{hdef}.  We deduce from this and \eqref{Hjdef} that
\begin{align*}
\begin{split}
 H_j(p)+\overline{H_j(p)}=& 2\Re h(p)+O(\frac {\log p}{\log X^{\alpha_j}}), \\ 
 H_j(\varpi)\overline{H_j(\overline{\varpi}})+H_j(\overline{\varpi})\overline{H_j(\varpi)}+H_j(\varpi)H_j(\overline{\varpi})+\overline{H_j(\varpi)H_j(\overline{\varpi})})=&  2|h(\sqrt{p})|^2+2\Re h(\sqrt{p})^2+O(\frac {\log p}{\log X^{\alpha_j}}), \quad p=\varpi\overline{\varpi}. 
\end{split}
\end{align*}    

 The above implies that 
\begin{align*}
\begin{split}
 & \sum_{\substack{n_i\overline{n'_i} \in \mn\\ \varpi | n_in'_i \Rightarrow X^{\alpha_{i-1}}<N(\varpi) \leq X^{\alpha_i}}} \frac{H_j(n_i)\overline{H_j(n'_i)}}{\sqrt{N(n_in'_i)}w(n_i)w(\overline{n'_i})} \\
 =& \prod_{\substack{X^{\alpha_{i-1}/2}<p \leq X^{\alpha_i/2} \\ p \ \text{inert} }}\Big  (1+\frac {2\Re h(p)}{p}+O\Big( \frac 1{p^{2}}+\frac {\log p}{p\log X^{\alpha_j}} \Big) \Big ) \\
 & \times \prod_{\substack{X^{\alpha_{i-1}/2}<p \leq X^{\alpha_i/2} \\ p=\varpi \overline{\varpi} }}\Big  (1+\frac {2|h(\sqrt{p})|^2+2\Re h(\sqrt{p})^2}{p}+O\Big( \frac 1{p^{2}}+\frac {\log p}{p\log X^{\alpha_j}}  \Big) \Big ) \\
 =& \prod_{\substack{X^{\alpha_{i-1}/2}<p \leq X^{\alpha_i/2} }}\Big  (1+\frac {2\Re h(p)}{p}\cdot\frac {1-\chi^{(D_K)}(p)}{2}+\frac {2|h(\sqrt{p})|^2+2\Re h(\sqrt{p})^2}{p}\cdot \frac {1+\chi^{(D_K)}(p)}{2}+O\Big( \frac {\log p}{p\log X^{\alpha_j}}+\frac 1{p^{2}}  \Big) \Big ),
\end{split}
\end{align*}
  where the last equality above follows from \eqref{splitcriteria}. 

  Applying the well-known estimation that $1+x \leq e^x$ for any real $x$, we see that
\begin{align}
\label{sumHjmain}
\begin{split}
 & \sum_{\substack{n_i\overline{n'_i} \in \mn\\ \varpi | n_in'_i \Rightarrow X^{\alpha_{i-1}}<N(\varpi) \leq X^{\alpha_i}}} \frac{H_j(n_i)\overline{H_j(n'_i)}}{\sqrt{N(n_in'_i)}w(n_i)w(\overline{n'_i})} \\
 \leq &  \exp \Big  (\sum_{\substack{X^{\alpha_{i-1}/2}<p \leq X^{\alpha_i/2} }} \Big ( 1+\frac {2\Re h(p)}{p}\cdot\frac {1-\chi^{(D_K)}(p)}{2}+\frac {2|h(\sqrt{p})|^2+2\Re h(\sqrt{p})^2}{p}\cdot \frac {1+\chi^{(D_K)}(p)}{2}\Big )+O\Big( \frac {\log p}{p\log X^{\alpha_j}}+\frac 1{p^{2}}  \Big)  \Big ).
\end{split}
\end{align}

  Similarly, we have
\begin{align}
\label{sumsqureierror1}
\begin{split}
 &\sum_{\substack{n_i\overline{n'_i} \in \mn\\ \varpi | n_in'_i \Rightarrow X^{\alpha_{i-1}}<N(\varpi) \leq X^{\alpha_i}}} \frac{2^{\Omega_K(n_i) -e^2\alpha^{-3/4}_i}|H_j(n_i)\overline{H_j(n'_i)}|}{\sqrt{N(n_in'_i)}w(n_i)w(\overline{n'_i})} \\
\leq & 2^{-\alpha^{-3/4}_i}\exp \Big  (\sum_{\substack{X^{\alpha_{i-1}/2}<p \leq X^{\alpha_i/2} }} \Big ( 1+\frac {2\Re h(p)}{p}\cdot\frac {1-\chi^{(D_K)}(p)}{2}+\frac {2|h(\sqrt{p})|^2+2\Re h(\sqrt{p})^2}{p}\cdot \frac {1+\chi^{(D_K)}(p)}{2}\Big ) \\
& +O\Big( \frac {\log p}{p\log X^{\alpha_j}}+\frac 1{p^{2}}  \Big)  \Big ), 
\end{split}
\end{align}      
  where we obtain the last estimation above by observing that it follows from \eqref{merten} that
\begin{align*}
   \sum_{X^{\alpha_{i-1}/2}<p \leq X^{\alpha_i/2} }\frac 1{p} \leq \log \Big (\frac {\alpha_i}{\alpha_{i-1}}\Big )+O(1)=O(1).
\end{align*}

   We deduce from \eqref{SrbounddecompH}--\eqref{sumsqureierror1} that
\begin{align}
\label{Srbounddecompexp}
\begin{split}
 &\sum_{\substack{n_i\overline{n'_i} \in \mn}} \frac{H_j(n_i)\overline{H_j(n'_i)}}{\sqrt{N(n_in'_i)}w(n_i)w(\overline{n'_i})}  b_i(n_i)b_i(\overline{n'_i}) \\
 \ll &
 (1+O(2^{-\alpha^{-3/4}_i})) \\
 & \times \exp \Big  (\sum_{\substack{X^{\alpha_{i-1}/2}<p \leq X^{\alpha_i/2} }} \Big ( 1+\frac {2\Re h(p)}{p}\cdot\frac {1-\chi^{(D_K)}(p)}{2}+\frac {2|h(\sqrt{p})|^2+2\Re h(\sqrt{p})^2}{p}\cdot \frac {1+\chi^{(D_K)}(p)}{2}\Big )+O\Big( \frac {\log p}{p\log X^{\alpha_j}}+\frac 1{p^{2}}  \Big)  \Big ).
\end{split}
\end{align}
  
   We conclude from \eqref{Srbounddecomp}--\eqref{SrbounddecompM} and \eqref{Srbounddecompexp} that
\begin{align}
\label{prodEestmain}
\begin{split}
  S_r
\leq &  2^{- u\lceil 2^{u/2}\rceil/2 }e^{-10^3(A+1)a/(2\alpha_{j+1})}X\prod^j_{i=1}\Big (1+O(2^{-\alpha^{-3/4}_i})\Big ) \\
& \times \exp \Big (\sum_{p \leq X^{\alpha_j}/2  }\Big (\frac {2\Re h(p)}{p}\cdot\frac {1-\chi^{(D_K)}(p)}{2}+\frac {2|h(\sqrt{p})|^2+2\Re h(\sqrt{p})^2}{p}\cdot \frac {1+\chi^{(D_K)}(p)}{2}\Big )
 +O\Big(\sum_{p \leq X^{\alpha_j}/2}\Big (\frac {\log p}{p\log X^{\alpha_j}}+\frac 1{p^{2}}\Big ) \Big )\Big ) \\
\ll & 2^{- u\lceil 2^{u/2}\rceil/2 }e^{-10^3(A+1)a/(2\alpha_{j+1})}X \\
& \times \exp \Big (\sum_{p \leq X^{\alpha_j}/2  }\Big ( \frac {2\Re h(p)}{p}\cdot\frac {1-\chi^{(D_K)}(p)}{2}+\frac {2|h(\sqrt{p})|^2+2\Re h(\sqrt{p})^2}{p}\cdot \frac {1+\chi^{(D_K)}(p)}{2}\Big )
 +O\Big(\sum_{p \leq X^{\alpha_j}/2}\Big (\frac {\log p}{p\log X^{\alpha_j}}+\frac 1{p^{2}}\Big ) \Big )\Big ) .
\end{split}
\end{align}

    Note that 
\begin{align}
\label{sumpsqure}
  \sum_{p \leq X^{\alpha_j}/2}\frac 1{p^{2}} \ll 1.
\end{align}
  Also, by \eqref{mertenpartialsummation},
\begin{align}
\label{sumsplogp}
  \sum_{p \leq X^{\alpha_j}/2}\frac {\log p}{p\log X^{\alpha_j}} \ll 1.
\end{align}

 We deduce from \eqref{prodEestmain}--\eqref{sumsplogp} that
\begin{align*}
  S_r \ll & 2^{- u\lceil 2^{u/2}\rceil/2 }e^{-10^3(A+1)a/(2\alpha_{j+1})}X \exp \Big (\sum_{p \leq X^{\alpha_j}/2  }\Big ( \frac {2\Re h(p)}{p}\cdot\frac {1-\chi^{(D_K)}(p)}{2}+\frac {2|h(\sqrt{p})|^2+2\Re h(\sqrt{p})^2}{p}\cdot \frac {1+\chi^{(D_K)}(p)}{2}\Big ) \Big ) .
\end{align*}
  
  Inserting the above into \eqref{L2k}, we see that
\begin{align}
\label{L2k1}
\begin{split}
  & \sum_{l \in \mathcal{S}(j)\bigcap \mathcal{P}(u)} \big|L\big(\tfrac12+it_1, \lambda_{M_Kl} \big) \big|^{a_1} \cdots \big| L\big(\tfrac12+it_{k}, \lambda_{M_Kl} \big) \big|^{a_{k}}  \\
 \ll & e^{a2^{u/2+4}}\exp \left( \frac {2(A+1)a}{\alpha_j} \right)(\mathcal{J}-j)2^{- u\lceil 2^{u/2}\rceil/2 }e^{-10^3(A+1)a/(2\alpha_{j+1})}X \\
 & \times \exp \Big (\sum_{p \leq X^{\alpha_j}/2  }\Big ( \frac {2\Re h(p)}{p}\cdot\frac {1-\chi^{(D_K)}(p)}{2}+\frac {2|h(\sqrt{p})|^2+2\Re h(\sqrt{p})^2}{p}\cdot \frac {1+\chi^{(D_K)}(p)}{2}\Big )\Big ) \\
 \ll & e^{a2^{u/2+4}}\exp \left( -\frac {10^2(A+1)a}{\alpha_j} \right)2^{- u\lceil 2^{u/2}\rceil/2 }X \\
 & \times \exp \Big (\sum_{p \leq X^{\alpha_j}/2  }\Big (\frac {2\Re h(p)}{p}\cdot\frac {1-\chi^{(D_K)}(p)}{2}+\frac {2|h(\sqrt{p})|^2+2\Re h(\sqrt{p})^2}{p}\cdot \frac {1+\chi^{(D_K)}(p)}{2}\Big )\Big ), 
\end{split}
\end{align}
 where the last estimation above follows from the observation that
\begin{align*}
\begin{split}
 \mathcal{J}-j \leq  \frac{\log(1/\alpha_{j})}{\log 20}, \ j \geq 1. 
\end{split}
\end{align*} 
  
   Observe that we have $2|h(\sqrt{p})|^2+2\Re h(\sqrt{p})^2 \geq 0$, $1+\chi^{(D_K)} \geq 0$, so that
\begin{align}
\label{sumph2}
\begin{split}
  \sum_{p \leq X^{\alpha_j}/2  }\frac {2|h(\sqrt{p})|^2+2\Re h(\sqrt{p})^2}{p}\cdot \frac {1+\chi^{(D_K)}(p)}{2} \leq \sum_{p \leq X }\frac {2|h(\sqrt{p})|^2+2\Re h(\sqrt{p})^2}{p}\cdot \frac {1+\chi^{(D_K)}(p)}{2}
\end{split}
\end{align}   
 
   Note also that $|2\Re h(p)| \leq a$ so that by \eqref{merten}, we have
\begin{align}
\label{sumph1}
\begin{split}
  \sum_{p \leq X^{\alpha_j}/2  }\frac {2\Re h(p)}{p}\cdot\frac {1-\chi^{(D_K)}(p)}{2} \leq & \sum_{p \leq X }\frac {2|h(\sqrt{p})|^2+2\Re h(\sqrt{p})^2}{p}\cdot \frac {1+\chi^{(D_K)}(p)}{2}+a\log \frac {1}{\alpha_j}+O(1) \\
  \leq & \sum_{p \leq X }\frac {2|h(\sqrt{p})|^2+2\Re h(\sqrt{p})^2}{p}\cdot \frac {1+\chi^{(D_K)}(p)}{2}+\frac {a}{\alpha_j}+O(1)
\end{split}
\end{align}   

   We apply \eqref{sumph2} and \eqref{sumph1} in \eqref{L2k1} to see that 
\begin{align*}
\begin{split}
  & \sum_{l \in \mathcal{S}(j)\bigcap \mathcal{P}(u)} \big|L\big(\tfrac12+it_1, \lambda_{M_Kl} \big) \big|^{a_1} \cdots \big| L\big(\tfrac12+it_{k}, \lambda_{M_Kl} \big) \big|^{a_{k}}  \\
 \ll & e^{a2^{u/2+4}}\exp \left( -\frac {10(A+1)a}{\alpha_j} \right)2^{- u\lceil 2^{u/2}\rceil/2 }X \\
 & \times \exp \Big (\sum_{p \leq X }\Big (\frac {2\Re h(p)}{p}\cdot\frac {1-\chi^{(D_K)}(p)}{2}+\frac {2|h(\sqrt{p})|^2+2\Re h(\sqrt{p})^2}{p}\cdot \frac {1+\chi^{(D_K)}(p)}{2}\Big )\Big ). 
\end{split}
\end{align*}

   As the sum of the right side expression over $m$ and $j$ converges, we derive from the above that
\begin{align}
\label{sumovermjsimplified}
\begin{split}
  & \sum_{u=0}^{\log \log X/2}\sum_{j=0}^{\mathcal{J}}\sum_{l \in \mathcal{S}(j)\bigcap \mathcal{P}(u)} \big| L\big( \tfrac12+it_1, \lambda_{M_Kl} \big ) \big|^{a_1} \cdots \big| L\big( \tfrac12+it_{k}, \lambda_{M_Kl}  \big) \big|^{a_{k}} \\
   \ll & X\exp \Big (\sum_{p \leq X }\Big (\frac {2\Re h(p)}{p}\cdot\frac {1-\chi^{(D_K)}(p)}{2}+\frac {2|h(\sqrt{p})|^2+2\Re h(\sqrt{p})^2}{p}\cdot \frac {1+\chi^{(D_K)}(p)}{2}\Big )\Big ) .
\end{split}
\end{align}

   Note that
\begin{align}
\label{hexp}
\begin{split}
& \sum_{p \leq X }\Big (\frac {2\Re h(p)}{2p}+\frac {2|h(\sqrt{p})|^2+2\Re h(\sqrt{p})^2}{2p}\Big ) \\
= & \sum_{p \leq X }\frac {1}{2p}\Big ( \sum^{k}_{j=1}a_j\cos(2t_j\log p) + 
\Big(\sum^{k}_{j=1}a_j\cos(t_j\log p) \Big)^2+\frac 12\Big( \sum^{k}_{j=1}a_j\sin(t_j\log p) \Big)^2\Big ) \\
=& \sum_{p \leq X }\frac {1}{2p}\Big (\frac 12\sum^{k}_{j=1}a^2_j+\frac 12\sum^{k}_{j=1}a^2_j\cos^2(t_j\log p)+2\sum_{1 \leq i<j \leq k}a_ia_j\cos(t_i\log p)\cos(t_j\log p)\\
& +\sum_{1 \leq i<j \leq k}a_ia_j\sin(t_i\log p)\sin(t_j\log p) +\sum^{k}_{j=1}a_j\cos(2t_j\log p)\Big ) \\
=& \sum_{p \leq X }\frac {1}{2p}\Big (\frac 34\sum^{k}_{j=1}a^2_j+\sum_{1 \leq i<j \leq k}a_ia_j\big(\frac 12\cos((t_i+t_j)\log p)+\frac 32\cos((t_i-t_j)\log p)\big )+\sum^{k}_{j=1} \Big( \frac {a^2_j}{4}+a_j \Big) \cos(2t_j\log p)\Big ),
\end{split}
\end{align}
upon using the relations
\begin{align*}
\begin{split}
 \cos (\alpha)\cos (\beta)=\frac 12 (\cos (\alpha+\beta)+\cos(\alpha-\beta)), \ \sin (\alpha)\sin (\beta)=-\frac 12 (\cos (\alpha+\beta)-\cos(\alpha-\beta)). 
\end{split}
\end{align*}

  Similarly, we have
\begin{align}
\label{hexp1}
\begin{split}
& \sum_{p \leq X }\Big (-\frac {2\Re h(p)}{2p}\chi^{(D_K)}(p)+\frac {2|h(\sqrt{p})|^2+2\Re h(\sqrt{p})^2}{2p}\chi^{(D_K)}(p)\Big ) \\
= & \sum_{p \leq X }\frac {\chi^{(D_K)}(p)}{2p}\Big ( -\sum^{k}_{j=1}a_j\cos(2t_j\log p) + 
\Big(\sum^{k}_{j=1}a_j\cos(t_j\log p) \Big)^2+\frac 12\Big( \sum^{k}_{j=1}a_j\sin(t_j\log p) \Big)^2\Big ) \\
=& \sum_{p \leq X }\frac {\chi^{(D_K)}(p)}{2p}\Big (\frac 34\sum^{k}_{j=1}a^2_j+\sum_{1 \leq i<j \leq k}a_ia_j\big(\frac 12\cos((t_i+t_j)\log p)+\frac 32\cos((t_i-t_j)\log p)\big )+\sum^{k}_{j=1} \Big( \frac {a^2_j}{4}-a_j \Big) \cos(2t_j\log p)\Big ).
\end{split}
\end{align}

   We now apply \eqref{mertenpartialsummation} and \eqref{mertenstype} to evaluate the last expression in \eqref{hexp}. We also apply \eqref{mertenpartialsummation}
     and \eqref{mertenstypesympower} to evaluate the last expression in \eqref{hexp1}. We then substitute the result into \eqref{sumovermjsimplified}.  
     This allows us to derive the desired estimation in \eqref{sumovermj} and hence completes the proof of Theorem \ref{t1}. 
   
\section{Proof of Theorem \ref{quadraticmean}}
\label{sec: mainthmquad}

   We begin with an auxiliary  result.
\begin{proposition}
\label{t3prop}
 With the notation as above and assuming the truth of GRH, we have for any fixed integer $k \geq 1$ and any real numbers $2m \geq 4k+4$, $10 \leq B=X^{O(1)}$,
\begin{align}
\label{finiteintestquad}
\begin{split}
  \sum_{1 \leq l \leq X} & \bigg(\int\limits_{0}^{B}|L(\tfrac{1}{2}+it, \lambda_{M_Kl})| \dif t\bigg)^{2m} \\
 \ll & X\big( (\log X)^{2m^2-m+3k+1}B^{4k}(\log \log B)^{O(1)}+(\log X)^{3(m-2k)^2/2+8k^2+k+1}B^{2m}(\log \log B)^{O(1)}(\log\log X)^{O(1)} \\
 & \hspace*{2cm} +(\log X)^{2(m-2k)^2+2m+8k^2-8k}B^{2m-1}(\log \log B)^{O(1)}(\log\log X)^{O(1)}\big).
\end{split}
\end{align}
\end{proposition}
\begin{proof}
  We first apply H\"older's inequality to see that
\begin{align*}
    \bigg(\int\limits_{0}^{B} |L(\tfrac{1}{2}+ it, \lambda_{M_Kl})| \dif t\bigg)^{2m}=& \bigg(\int\limits_{0}^{B} |L(\tfrac{1}{2}+ it, \lambda_{M_Kl})| \dif t\bigg)^{4k}\bigg(\int\limits_{0}^{B} |L(\tfrac{1}{2}+ it, \lambda_{M_Kl})| \dif t\bigg)^{2m-4k} \\
      \ll & B^{3k} \bigg(\int\limits_{0}^{B} |L(\tfrac{1}{2}+ it, \lambda_{M_Kl})|^4 \dif t\bigg)^{k} \bigg(\int\limits_{0}^{B} |L(\tfrac 12+iu, \lambda_{M_Kl})| \dif u \bigg)^{2m-4k}.
\end{align*}  
   It follows from the above and symmetry that we have 
\begin{align}
\label{Lintdecomp}
    \bigg(\int\limits_{0}^{B} |L(\tfrac{1}{2}+ it, \lambda_{M_Kl})| \dif t\bigg)^{2m}
      \ll B^{3k}\int\limits_{[0,B]^k}\prod_{a=1}^k|L(\tfrac 12+ it_a, \lambda_{M_Kl})|^4 \bigg(\int_{\mathcal{D} }|L(\tfrac 12+iu, \lambda_{M_Kl})| \dif u \bigg)^{2m-4k} \dif \mathbf{t},
\end{align}
where $\mathcal{D}=\mathcal{D}(t_1,\ldots,t_k)=\{ u\in [0,B]:|t_1-u|\leq |t_2-u|\leq \ldots \leq |t_k-u| \}$. 

 We set $\mathcal{B}_1=\big[-\frac{1}{\log X},\frac{1}{\log X}\big]$, $\mathcal{B}_j=\big[-\frac{e^{j-1}}{\log X}, -\frac{e^{j-2}}{\log X}\big]
  \cup \big[\frac{e^{j-2}}{\log X}, \frac{e^{j-1}}{\log X}\big]$ for $2\leq j< \lfloor \log \log X\rfloor+10 =: W$ and 
  $\mathcal{B}_W=[-B,B]\setminus \bigcup_{1\leq j<W} \mathcal{B}_j$. 

As we have $\mathcal{D}\subset [0,B] \subset t_1+[-B,B]\subset \bigcup_{1\leq j\leq W} ( t_1+\mathcal{B}_j)$ for any $t_1\in [0,B]$, we see that, by
setting $\mathcal{A}_j=\mathcal{B}_j\cap (-t_1+\mathcal{D})$, the sets $(t_1+\mathcal{A}_j)_{1\leq j\leq W}$ form a partition of $\mathcal{D}$.  We then 
apply Hölder's inequality twice to obtain that
\begin{align}
\label{LintoverD}
\begin{split}
    \bigg(\int\limits_{\mathcal{D}} & |L(\tfrac{1}{2} + iu, \lambda_{M_Kl})| \dif u\bigg)^{2m-4k}  \leq \bigg( \sum_{1\leq j\leq W} \frac{1}{j}\cdot  j \int\limits_{t_1+\mathcal{A}_j} |L( \tfrac{1}{2}+iu, \lambda_{M_Kl})| \dif u  \bigg)^{2m-4k} \\
     \leq & \bigg(\sum_{1\leq j\leq W} j^{2m-4k} \bigg( \int\limits_{t_1+\mathcal{A}_j} \big|L( \tfrac{1}{2}+iu, \lambda_{M_Kl})\big| \dif u  \bigg)^{2m-4k}\bigg)
     \bigg(\sum_{1\leq j\leq W } j^{-(2m-4k)/(2m-4k-1)} \bigg)^{2m-4k-1} \\
     \ll & \sum_{1\leq j\leq W} j^{2m-4k} \bigg( \int\limits_{t_1+\mathcal{A}_j} |L( \tfrac{1}{2}+iu, \lambda_{M_Kl})| \dif u \bigg)^{2m-4k} \\
     \leq & \sum_{1\leq j\leq W} j^{2m-4k} |\mathcal{B}_j|^{2m-4k-1} \int\limits_{t_1+\mathcal{A}_j} |L(\tfrac12+iu, \lambda_{M_Kl})|^{2m-4k} \dif u.
\end{split}
\end{align}
  For $\mathbf{t}=(t_1,\ldots,t_k)$,  write
$$L(\mathbf{t},u) =\sum_{1 \leq l \leq X} 
\prod_{a=1}^k|L( \tfrac{1}{2}+it_a, \lambda_{M_Kl})|^4 \cdot |L( \tfrac{1}{2}+iu, \lambda_{M_Kl})|^{2m-4k}.$$
 We then deduce from \eqref{Lintdecomp} and \eqref{LintoverD} that
\begin{align}
\label{Lintest}
\begin{split}
    \sum_{1 \leq l \leq X} \bigg(\int\limits_{0}^{B}|L(\tfrac 12+it,\lambda_{M_Kl})|\dif t \bigg)^{2m}\ll &
    B^{3k}\sum_{1\leq l_0\leq W} l_0^{2m-4k} |\mathcal{B}_{l_0}|^{2m-4k-1} \int\limits_{[0,B]^k}\int\limits_{t_1+\mathcal{A}_{l_0}} L(\mathbf{t},u) \dif u \dif \mathbf{t}  \\
     \ll & B^{3k}\sum_{1\leq l_0, l_1, \ldots l_{k-1}\leq W} l_0^{2m-4k} |\mathcal{B}_{l_0}|^{2m-4k-1} \int\limits_{\mathcal{C}_{l_0,l_1, \cdots, l_{k-1}}} L(\mathbf{t},u) \dif u \dif \mathbf{t},
\end{split}
\end{align}
where
$$\mathcal{C}_{l_0,l_1, \cdots, l_{k-1}}=\{(t_1,\ldots,t_k,u)\in [0,B]^{k+1}: u\in t_1+ \mathcal{A}_{l_0},\, |t_{i+1}-u|-|t_i-u|\in \mathcal{B}_{l_i}, \ 1 \leq i \leq k-1\}.$$
  Now we devide the last summation of \eqref{Lintest} into two cases depending on the size of $l_0$. 

\textbf{Case 1:} $l_0<W$. Note that the volume of $\mathcal{C}_{l_0,l_1, \cdots, l_{k-1}}$ is $\ll  B^k e^{l_0+l_1+\cdots+l_{k-1}} (\log X)^{-k}$. 
By the definition of $\mathcal{C}_{l_0,l_1, \cdots, l_{k-1}}$ and the observation $t_i$, $u \geq 0$, $1\leq i \leq k$ 
that we have $e^{l_0}/\log X \ll |t_1-u|\ll |t_1+u|\ll B =X^{O(1)}$. This implies that $g_1(|t_1\pm u|)\ll \log X \cdot \log \log B/e^{l_0}$, where 
the function $g_1(x)$ is defined in \eqref{g3Def}.  It also follows from the definition of $\mathcal{A}_j$ that $|t_2-u|\geq |t_1-u|$, so
 that $B \gg |t_2+u| \gg |t_2-u|= |t_1-u|+(|t_2-u|-|t_1-u|)\gg  e^{l_0}/\log X + e^{l_1}/\log X$. This implies that $g_1(|t_2\pm u|)\ll \log X \cdot \log \log B/e^{\max(l_0,l_1) }$.
Similarly, we have  $g_1(|t_i\pm u|)\ll \log X \cdot \log \log B /e^{\max(l_0,l_1,\ldots, l_{i-1}) }$ for any $1 \leq i \leq k$.
Moreover, we have $\sum^{j-1}_{s=i}(|t_{s+1}-u|-|t_s-u|) \leq |t_j-t_i| \leq |t_j+t_i|$ for any $1 \leq i < j \leq k$, which implies that $g_1(|t_{j}\pm t_i|)\ll \log X \cdot \log \log B/e^{\max(l_i,\ldots, l_{j-1} ) }$.
We now bound $g_1(|2t_i|), 1\leq i \leq k$, $g_1(|2u|)$ trivially by $\log X$ and the function $g_2$ given in \eqref{g4Def} by $(\log \log B)^{O(1)}$ 
to derive from Corollary \ref{cor1quad} that when $2m \geq 4k+4$, we have for $(t_1,\ldots,t_k,u)\in \mathcal{C}_{l_0,l_1, \cdots, l_{k-1}}$,
\begin{align*}
     & L(\mathbf{t},u) \\
     \ll & X(\log X)^{((2m-4k)^2+4^2k)/2+((2m-4k)+4k)/2} \\
     &  \hspace*{2cm} \times \bigg(\prod^{k-1}_{i=0}\frac{\log X}{e^{ \max(l_0,l_1,\ldots, l_{i}) }} \bigg)^{4(2m-4k)}
   \bigg(\prod^{k-1}_{i=1} \prod^{k}_{j=i+1}\frac{\log X}{e^{\max(l_i,\ldots, l_{j-1} ) } } \bigg)^{16} (\log \log B)^{O(1)} \\
     = & X(\log X)^{m(2m+1)} \exp\Big( -4(2m-4k)\sum^{k-1}_{i=0}\max(l_0,l_1,\ldots, l_{i})-16\sum^{k-1}_{i=1} \sum^{k}_{j=i+1}\max(l_i,\ldots, l_{j-1} )\Big) (\log \log B)^{O(1)} .
\end{align*}

As $|\mathcal{B}_{l_0}|\ll e^{l_0}/\log X$, we see that
\begin{align}
\label{firstcase}
\begin{split}
       &  \sum_{\substack{1\leq l_0<W \\ 1\leq l_1, \ldots l_{k-1}\leq W}}  l_0^{2m-4k} |\mathcal{B}_{l_0}|^{2m-4k-1} \int\limits_{\mathcal{C}_{l_0,l_1, \cdots, l_{k-1}}} L(\mathbf{t},u) \dif u \dif \mathbf{t} \\
    \ll & X(\log X)^{2m^2-m+3k+1}B^k(\log \log B)^{O(1)}  \cdot \\
    &  \times \sum_{\substack{1\leq l_0<W \\ 1\leq l_1, \ldots l_{k-1}\leq W}}  l_0^{2m-4k}\exp\Big( (2m-4k-1)l_0+\sum^{k-1}_{i=0}l_i-4(2m-4k)\sum^{k-1}_{i=0}\max(l_0,l_1,\ldots, l_{i})-16\sum^{k-1}_{i=1} \sum^{k}_{j=i+1}\max(l_i,\ldots, l_{j-1} )\Big) \\
    = & X(\log X)^{2m^2-m+3k+1}B^k(\log \log B)^{O(1)}  \cdot \\
    &  \times \sum_{\substack{1\leq l_0<W \\ 1\leq l_1, \ldots l_{k-1}\leq W}}  l_0^{2m-4k}\exp\Big( -3(2m-4k)l_0-15\sum^{k-1}_{i=1}l_i-4(2m-4k)\sum^{k-1}_{i=1}\max(l_0,l_1,\ldots, l_{i})-\sum^{k-1}_{i=1} \sum^{k}_{j=i+2}\max(l_i,\ldots, l_{j-1} )\Big) \\
    \ll &   X(\log X)^{2m^2-m+3k+1}B^k(\log \log B)^{O(1)}.
\end{split}
\end{align}

\textbf{Case 2:} $l_0=W$.  Observe that $g_1(|t_i\pm u|)\ll \log\log B$ for each $1 \leq i \leq k$. Also, similar to Case 1, 
we have $g_1(|t_j\pm t_i|) \ll \log X/e^{\max(l_i,\ldots, l_{j-1} ) }$ for $1 \leq i \leq k$. When $|u| \geq 5$, we have $g_1(|2u|) \ll \log\log B$. 
We now bound the functions $g_2$ by $(\log \log B)^{O(1)}$ and $g_1(|2t_i|)$ by $\log X$ for $1 \leq i \leq k$ to derive from Corollary \ref{cor1quad} that
\begin{align*}
     L(\mathbf{t},u) \ll & X(\log X)^{3(2m-4k)^2/8+4^2k/2+4k/2}
     \bigg(\prod^{k-1}_{i=1} \prod^{k}_{j=i+1}\frac{\log X}{e^{\max(l_i,\ldots, l_{j-1} ) } } \bigg)^{16} (\log \log B)^{O(1)} \\
     = & X(\log X)^{3(m-2k)^2/2+8k^2+2k} \exp\Big( -16\sum^{k-1}_{i=1} \sum^{k}_{j=i+1}\max(l_i,\ldots, l_{j-1} )\Big) (\log \log B)^{O(1)}.
\end{align*}

  If $|u| \leq 5$, then as $|t_1-u| \leq |t_i-u| \leq |t_i|+|u|$ for any $1 \leq i \leq k$, we see that $|t_i| \geq 5$, $g_1(|2t_i|) \ll \log\log B$.  
  Bounding $g_1(|2u|)$ by $\log X$ and the function $g_2$ by $(\log \log B)^{O(1)}$ again, we derive from Corollary \ref{cor1quad} that
\begin{align*}
 L(\mathbf{t},u) \ll & X(\log X)^{3((2m-4k)^2+4^2k)/8+(2m-4k)^2/8+(2m-4k)/2}
     \bigg(\prod^{k-1}_{i=1} \prod^{k}_{j=i+1}\frac{\log X}{e^{\max(l_i,\ldots, l_{j-1} ) } } \bigg)^{16} (\log \log B)^{O(1)} \\
     = & X(\log X)^{2(m-2k)^2+2m+8k^2-8k} \exp\Big( -16\sum^{k-1}_{i=1} \sum^{k}_{j=i+1}\max(l_i,\ldots, l_{j-1} )\Big) (\log \log B)^{O(1)} .
\end{align*}
  Observe that the volume of the region $\mathcal{C}_{W,l_1, \cdots, l_{k-1}}$ is $\ll B^{k+1} e^{l_1+\cdots+l_{k-1}} (\log X)^{-k+1}$.  
  As $|\mathcal{B}_W|\ll B$, we derive from the above that
\begin{align}
\label{secondcase}
\begin{split}
     \sum_{\substack{1\leq l_1, \ldots l_{k-1}\leq W}} & W^{2m-4k} |\mathcal{B}_{W}|^{2m-4k-1} \int\limits_{\substack{\mathcal{C}_{K,l_1, \cdots, l_{k-1}} \\ |u| \geq 5}} L(\mathbf{t},u) \dif u \dif \mathbf{t}  \\
     \ll & X(\log X)^{3(m-2k)^2/2+8k^2+k+1}B^{2m-3k}(\log \log B)^{O(1)}(\log\log X)^{O(1)} \\
     & \hspace*{2cm} \times \sum_{\substack{1\leq l_1, \ldots l_{k-1}\leq W}} \exp\Big( \sum^{k-1}_{i=1}l_i-16\sum^{k-1}_{i=1} \sum^{k}_{j=i+1}\max(l_i,\ldots, l_{j-1} )\Big) \\
     \ll &   X(\log X)^{3(m-2k)^2/2+8k^2+k+1}B^{2m-3k}(\log \log B)^{O(1)}(\log\log X)^{O(1)}.
\end{split}
\end{align}

   Notice that when $|u| \leq 5$, the volume of the region $\mathcal{C}_{K,l_1, \cdots, l_{k-1}}$ is $\ll B^{k}$.  Thus
\begin{align}
\label{secondcaseusmall}
\begin{split}
 \sum_{\substack{1\leq l_1, \ldots l_{k-1}\leq W}}  W^{2m-4k} & |\mathcal{B}_{W}|^{2m-4k-1} \int\limits_{\substack{\mathcal{C}_{K,l_1, \cdots, l_{k-1}} \\ |u| \leq 5}} L(\mathbf{t},u) \dif u \dif \mathbf{t}  \\
     \ll & X(\log X)^{2(m-2k)^2+2m+8k^2-8k}B^{2m-3k-1}(\log \log B)^{O(1)}(\log\log X)^{O(1)} \\
     & \hspace*{2cm} \times \sum_{\substack{1\leq l_1, \ldots l_{k-1}\leq W}} \exp\Big( -4\sum^{k-1}_{i=1} \sum^{k}_{j=i+1}\max(l_i,\ldots, l_{j-1} )\Big) \\
     \ll &   X(\log X)^{2(m-2k)^2+2m+8k^2-8k}B^{2m-3k-1}(\log \log B)^{O(1)}(\log\log X)^{O(1)}.
\end{split}
\end{align}
     
  It follows from \eqref{Lintest}--\eqref{secondcaseusmall} that the estimation given in \eqref{finiteintestquad} is valid. This completes the proof of the proposition.
\end{proof}

   Now, let $\Phi(t)$ be as defined in Section \ref{sec 1}.  Note that repeated integration by parts implies that the Mellin transform $\widehat{\Phi}$ of $\Phi$ 
   satisfies the bound that  for any integer $E \geq 1$ and $\Re(s) \geq 1/2$,
\begin{align}
\label{whatbound}
 \widehat{\Phi}(s)  \ll  X^{(E-1)\varepsilon_0}(1+|s|)^{-E}.
\end{align}

  We apply the Mellin inversion to see that
\begin{align}
\label{charintinitial1}
\begin{split}
 S_m(X,Y;K, \Phi) =&
\sum_{1 \leq l \leq X} \Big | \frac 1{2\pi i}\int\limits_{(2)}L(s, \lambda_{M_Kl}) Y^s\widehat{\Phi}(s) \dif s\Big |^{2m}.
\end{split}
\end{align}
  
  It follows from \cite[Corollary 5.20]{iwakow} that for every $1 \leq l \leq X$, $\Re(s) \geq 1/2$ and any $\varepsilon>0$, we have under GRH,
\begin{align}
\label{Lbound}
 L(s, \lambda_{M_Kl}) \ll |s|^{\varepsilon}.
\end{align}

Based on \eqref{whatbound} and \eqref{Lbound}, we may shift the line of integration in \eqref{charintinitial1} to $\Re(s)=1/2$ so that
\begin{align}
\label{charint}
\begin{split}
   S_m(X,Y;K, \Phi) =&
\sum_{1 \leq l \leq X} \Big | \frac 1{2\pi i}\int\limits_{(1/2)}L(s, \lambda_{M_Kl}) Y^s\widehat{\Phi}(s) \dif s\Big |^{2m}. 
\end{split}
\end{align} 	

 We split the integral in \eqref{charint} into two ranges, one over $|\Im(s)| \leq X^{4\varepsilon_0}$ and the other over $|\Im(s)| > X^{4\varepsilon_0}$ to see that
\begin{align*}
 S_m(X,Y;K, \Phi) \ll \sum_{1 \leq l \leq X} \Big | \int\limits_{\substack{ (1/2) \\ |\Im(s)| \leq X^{4\varepsilon_0}}}L(s, \lambda_{M_Kl})Y^s\widehat{\Phi}(s) \dif s\Big |^{2m}+
\sum_{1 \leq l \leq X} \Big | \int\limits_{\substack{ (1/2) \\ |\Im(s)| > X^{4\varepsilon_0}}}L(s, \lambda_{M_Kl})Y^s\widehat{\Phi}(s) \dif s\Big |^{2m}.
\end{align*} 	

   We apply H\"older's inequality with the observation that $m>2$ to see that
\begin{align}
\label{LintImslarge}
\begin{split}
 \sum_{1 \leq l \leq X} \Big | \int\limits_{\substack{ (1/2) \\ |\Im(s)| > X^{4\varepsilon_0}}}L(s, \lambda_{M_Kl})Y^s\widehat{\Phi}(s) \dif s\Big |^{2m}  \ll Y^m\Big (\int\limits_{\substack{ (1/2) \\ |\Im(s)| > X^{4\varepsilon_0}}}\Big | \widehat{\Phi}(s) \Big| |\dif s| \Big )^{2m-1}
 \int\limits_{\substack{ (1/2) \\ |\Im(s)| > X^{4\varepsilon_0}}}\sum_{1 \leq l \leq X} \Big |L(s, \lambda_{M_Kl})\Big|^{2m} \Big| \widehat{\Phi}(s)\Big | |\dif s|.
\end{split}
\end{align} 	

By \eqref{whatbound}, we have that
\begin{align}
\label{phiint}
\begin{split}
  \int\limits_{\substack {(1/2) \\ |\Im(s)| >  X^{4\varepsilon_0}}}\Big | \widehat{\Phi}(s) \Big| |\dif s|
 \ll & \int\limits_{\substack {(1/2) \\ |\Im(s)| >  X^{4\varepsilon_0}}}\frac  {X^{\varepsilon_0}}{1+|s|^2} |\dif s| \ll_{\varepsilon_0} \frac {X^{\varepsilon_0}}{X^{4\varepsilon_0}}.
\end{split}
\end{align}

 We also apply \eqref{Lbound} to see that
\begin{align}
\label{Lintslarge}
\begin{split}
  & \int\limits_{\substack {(1/2) \\ |\Im(s)| > X^{4\varepsilon_0}}}\sum_{1 \leq l \leq X}
  \Big | L(s, \lambda_{M_Kl})\Big |^{2m} \cdot \Big |\widehat{\Phi}(s)\Big | |\dif s| \ll  \int\limits_{\substack {(1/2) \\ |\Im(s)| >X^{4\varepsilon_0} }}\sum_{1 \leq l \leq X} |s|^{\varepsilon}\cdot\frac  {X^{\varepsilon_0}}{1+|s|^2}  |\dif s| \ll X^{1+\varepsilon_0}X^{-4\varepsilon_0(1-\varepsilon)+\varepsilon}.
\end{split}
\end{align}

  It follows from \eqref{charintinitial1}, \eqref{charint}--\eqref{Lintslarge} that
\begin{align}
\label{Ssimplified}
\begin{split}
 S_m(X,Y;K, \Phi) \ll &
    \sum_{1 \leq l \leq X} \Big | \int\limits_{\substack{ (1/2) \\ |\Im(s)| \leq X^{4\varepsilon_0}}}L(s, \lambda_{M_Kl})Y^s\widehat{\Phi}(s) \dif s\Big |^{2m}+
   O(XY^m) \\
   \ll & Y^m \sum_{1 \leq l \leq X}
   \Big | \int\limits_{\substack{ |t| \leq X^{4\varepsilon_0}}}\Big |L( \tfrac{1}{2}+it, \lambda_{M_Kl})\Big |\frac 1{1+|t|} \dif t\Big |^{2m}
   +O(XY^m). 
\end{split}
\end{align}
  
    We deduce from the above by applying H\"older's
inequality and symmetry that for $a=1-1/(2m)+\varepsilon$ with $\varepsilon>0$,
\begin{align*}
\begin{split}
 \Big | \int\limits_{ |t| \leq X^{4\varepsilon_0}}   \Big |L(\tfrac{1}{2}+it, \lambda_{M_Kl})|\frac {\dif t}{1+|t|} \Big |^{2m} 
 \ll & \Big |\int\limits_0^{X^{4\varepsilon_0}} \frac{|L(\tfrac{1}{2}+it,\lambda_{M_Kl})|}{t+1} \dif t\Big |^{2m} \\
  \leq & \bigg(\sum_{n\leq \log X+1} n^{-2am/(2m-1)} \bigg)^{2m-1}
    \sum_{n\leq  \log X+1} \bigg(n^a\int\limits_{e^{n-1}-1}^{e^{n}-1 } \frac{|L(\tfrac{1}{2}+it,\lambda_{M_Kl}) |}{t+1} \dif t\bigg)^{2m}   \\
 \ll  &  \sum_{n\leq  \log X+1} \frac{n^{2m-1+\varepsilon} }{e^{2nm} } \bigg( \int\limits_{e^{n-1}-1}^{e^{n}-1 } |L(\tfrac{1}{2}+it,\lambda_{M_Kl}) | \dif t \bigg)^{2m}.
\end{split}
\end{align*}

We now apply Proposition \ref{t3prop} to obtain that for any integer $k \geq 1$ and any real numbers $2m \geq 4k+4, \varepsilon>0$,
\begin{align*}
 & \sum_{n\leq  \log X+1}  \frac{n^{2m-1+\varepsilon} }{e^{4nm} }\sum_{1 \leq l \leq X} 
  \bigg( \int\limits_{e^{n-1}-1}^{e^{n}-1 } |L(\tfrac12+it,\lambda_{M_Kl}) | \dif t \bigg)^{2m}
    \\
     \ll & X\sum_{n\leq  \log X+1} \frac{n^{2m-1+\varepsilon} }{e^{2nm} } \\
    & \times
     \Big((\log X)^{2m^2-m+3k+1}e^{4kn}(\log 2n)^{O(1)}+(\log X)^{3(m-2k)^2/2+8k^2+k+1}(\log 2n)^{O(1)}(\log\log X)^{O(1)} e^{2mn} \\
    & \hspace*{2cm} +(\log X)^{2(m-2k)^2+2m+8k^2-8k}(\log 2n)^{O(1)}(\log\log X)^{O(1)} e^{(2m-1)n}\Big)  \ll X(\log X)^{E(m,k,\varepsilon)},
\end{align*}
  where 
\begin{equation}
\label{Edef}
 E(m,k,\varepsilon) = \max (2m^2-m+3k+1, 3(m-2k)^2/2+8k^2+k+1+2m+\varepsilon, 2(m-2k)^2+2m+8k^2-8k+\varepsilon) .
\end{equation}

 We set $k=1$ in \eqref{Edef} to deduce that the estimation given in \eqref{mainestimationmlarge} is valid. This completes the proof of Theorem \ref{quadraticmean}.

\vspace*{.5cm}

\noindent{\bf Acknowledgments.}  P. G. is supported in part by NSFC grant 12471003.
\bibliography{biblio}
\bibliographystyle{amsxport}

\end{document}